\documentclass[12pt]{article}
\usepackage{t1enc}
\usepackage{amsmath}
\usepackage{graphicx}
\usepackage{amsfonts}
\usepackage{amssymb}

\begin{document}

\title{Orbits of Families of Vector Fields\\on Subcartesian Spaces }
\author{J\k{e}drzej \'{S}niatycki\\Department of Mathematics and Statistics\\University of Calgary\\Calgary, Alberta, Canada\\e-mail: sniat@math.ucalgary.ca}
\maketitle
\begin{abstract}
Orbits of complete families of vector fields on a subcartesian space are shown
to be smooth manifolds. This allows a description of the structure of the
reduced phase space of a Hamiltonian system in terms of the reduced Poisson algebra.

Moreover, one can give a global description of smooth geometric structures on
a family of manifolds, which form a singular foliation of a subcartesian
space, in terms of objects defined on the corresponding family of vector
fields. Stratified spaces, Poisson spaces, and almost complex spaces are
discussed as examples.
\end{abstract}

\bigskip

\bigskip



\bigskip

\noindent To appear in Annales de l'Institut Fourier.

\newpage%

\noindent
\textbf{Contents}

\medskip

1. Introduction

2. Review of differential spaces

3. Properties of subcartesian spaces

4. Families of vector fields

5. Orbits and integral manifolds

6. Stratified spaces

7. Poisson reduction

8. Subcartesian Poisson spaces

9. Almost complex structures

\section{Introduction}

This work is motivated by the program of Poisson reduction of Hamiltonian
systems. Under the assumption that the action of the symmetry group $G$ on the
phase space $P$ of the system is proper, the orbit space $S=P/G$ is stratified
by orbit type, \cite{D-K}. For a Hamiltonian system, each stratum of $S$ is
singularly foliated by symplectic leaves, \cite{bates-lerman},
\cite{sjamaar-lerman}. The orbit space $S$ has a differential structure
$C^{\infty}(S)$ given by push-forwards to $S$ of $G$-invariant smooth
functions on $P,$ \cite{schwarz}, \cite{cushman-sniatycki}. Moreover,
$C^{\infty}(S)$ has the structure of a Poisson algebra, usually called the
reduced Poisson algebra. Following the approach initiated by Sjamaar and
Lerman, \cite{sjamaar-lerman}, we want to describe strata of the
stratification of $S$ as well as leaves of the singular foliation directly in
terms of the reduced Poisson algebra $C^{\infty}(S)$.

The essential property of a smooth stratified space needed here is the fact
that it is a smooth subcartesian space. The notion of a subcartesian space was
introduced by Aronszajn, \cite{aronszajn}, and subsequently developed by
Aronszajn and Szeptycki, \cite{aronszajn-szeptycki1},
\cite{aronszajn-szeptycki}, and by Marshall, \cite{marshall}, \cite{marshall1}%
. Related notions were independently introduced and studied by Spallek,
\cite{spallek1}, \cite{spallek2}.

A smooth subcartesian space is a diffferential space in the sense of
Sikor\-ski, \cite{sikorski1}, \cite{sikorski2}, \cite{sikorski}, that is
locally diffeomorphic to a subset of a Cartesian space $\mathbb{R}^{n}.$
Hence, we can use the differential space approach and study properties of a
subcartesian space in terms of its ring of globally defined smooth functions.

In this paper, we generalize to smooth subcartesian spaces the theorem of
Sussmann on orbits of families of vector fields on manifolds, \cite{sussmann},
and investigate its applications. In order to do this, we must first extend to
subcartesian spaces the results on the relationship between derivations and
local one-parameter local groups of diffeomorphisms of locally semi-algebraic
differential spaces obtained in \cite{sniatycki 2002}.

Let $S$ be a smooth subcartesian space, and $X:C^{\infty}(S)\rightarrow
C^{\infty}(S):h\mapsto X\cdot h$ be a derivation of $C^{\infty}(S)$. A curve
$c:I\rightarrow S,$ where $I$ is an interval in $\mathbb{R}$, is an integral
curve of $X$ if
\[
\frac{d}{dt}h(c(t))=(X\cdot h)(c(t))\quad\text{for all }\,\,h\in C^{\infty
}(S),\,\,t\in I.
\]
We show that, for every derivation $X$ of $C^{\infty}(S)$ and every $x\in S$,
there exists a unique maximal integral curve of $X$ passing through $x$.

We define a vector field on a smooth subcartesian space $S$ to be derivation
that generates a local one-parameter group of local diffeomorphisms of $S$.
Let $\mathcal{F}$ be a family of vector fields on $S$. An orbit of
$\mathcal{F}$ through a point $x\in S$ is the maximal set of points in $S$
which can be joined to $x$ by piecewise smooth integral curves of vector
fields in $\mathcal{F}$. In other words, $y\in S$ belongs to the orbit of
$\mathcal{F}$ through $x$ if there exists a positive integer $m$, vector
fields $X^{1},...,X^{m}\in\mathcal{F}$ and $(t_{1},...,t_{m})\in\mathbb{R}%
^{m}$ such that
\[
y=\left(  \varphi_{t_{m}}^{X^{m}}\raisebox{2pt}{$\scriptstyle\circ
\, $}...\raisebox{2pt}{$\scriptstyle\circ\, $}\varphi_{t_{1}}^{X^{1}}\right)
(x).
\]
We introduce the notion of a locally complete family of vector fields. A
family $\mathcal{F}$ is locally complete if, for every $X,Y\in\mathcal{F}$,
$t\in\mathbb{R},$ and $x\in S,$ for which the push-forward $\varphi_{t\ast
}^{X}Y(x)$ is defined, there exists an open neighbourhood $U$ of $x$ and a
vector field $Z\in\mathcal{F}$ such that the restriction of $\varphi_{t\ast
}^{X}Y$ of $U$ coincides with the restriction of $Z$ to $U$. In particular, a
family consisting of a single vector field $X$ on $S$ is locally complete
because, for every $t\in\mathbb{R},$ $\varphi_{t\ast}^{X}X$ coincides with the
restriction of $X$ to the domain of $\varphi_{t\ast}^{X}X$. Thus, the notion
of local completeness of a family of vector fields is unrelated to
completeness of vector fields constituting the family.

\begin{description}
\item [Main Theorem]Each orbit of a locally complete family of vector fields
on a smooth subcartesian space $S$ is a smooth manifold, and its inclusion
into $S$ is smooth.
\end{description}

\noindent We refer to the partition of $S$ by orbits of $\mathcal{F}$ as the
singular foliation of $S$ defined by $\mathcal{F}$, and to orbits of
$\mathcal{F}$ as leaves of the singular foliation. In the case when $S$ is a
smooth manifold, and $\mathcal{F}$ is a locally complete family of smooth
vector fields on $S$, orbits of $\mathcal{F}$ give rise to a singular
foliation of $S$ in the sense of Stefan, \cite{stefan}. Stefan's definition of
a singular foliation of a smooth manifold contains a condition of local
triviality, similar to a local triviality of a stratification (see Section 6).
Orbits of a locally complete family of vector fields on a subcartesian space
need not satisfy an obvious extension of Stefan's condition.

We show that the family $\mathcal{X}(S)$ of all vector fields on a
subcartesian space $S$ is locally complete. The singular foliation of $S$
defined by $\mathcal{X}(S)$ is minimal in the sense that, for every family
$\mathcal{F}$ of vector fields on $S$, orbits of $\mathcal{F}$ are contained
in orbits of $\mathcal{X}(S).$ In particular, the restriction of $\mathcal{F}
$ to each orbit $M$ of $\mathcal{X}(S)$ is a family $\mathcal{F}_{M}$ of
vector fields on $M$, and orbits of $\mathcal{F}$ contained in $M$ are orbits
of $\mathcal{F}_{M}.$

We show that smooth stratified spaces are subcartesian. All stratified spaces
considered here are assumed to be smooth. A stratified space $S$ is locally
trivial if it is locally diffeomorphic to the product of a stratified space
and a cone, \cite{pflaum}. We introduce the notion of a strongly stratified
vector field on a startified space, and prove that the family of all strongly
stratified vector fields on a locally trivial stratified space $S$ is locally
complete and that its orbits are strata of $S $. Hence, each stratum of $S$ is
contained in an orbit of the family $\mathcal{X}(S)$ of all vector fields on
$S$. Moreover, we show that if $S$ is a locally trivial stratified space then
orbits of the family $\mathcal{X}(S)$ of all vector fields on $S$ also give
rise to a stratification of $S.$ If the original stratification of $S$ is
minimal, then it coincides with the stratification by orbits of $\mathcal{X}%
(S)$. These results on stratified spaces are applied to describe singular
Poisson reduction of Hamiltonian systems.

We discuss also subcartesian Poisson spaces and almost complex spaces. A
combination of these two structures gives rise to a generalization to
subcartesian spaces of stratified K\"{a}hler spaces studied by Huebschmann
\cite{huebschmann}.

\section{Differential Spaces}

We begin with a review of elements of the theory of differential spaces,
\cite{sikorski}. Results stated here will be used in our study of vector
fields on subcartesian spaces.

A differential structure on a topological space $R$ is a family of functions
$C^{\infty}(R)$ satisfying the following conditions:

\begin{description}
\item [2.1.]The family
\[
\{f^{-1}((a,b))\mid f\in C^{\infty}(R),\,a,b\in\mathbb{R}\}
\]
is a sub-basis for the topology of $R.$

\item[2.2.] If $f_{1},...,f_{n}\in C^{\infty}(R)$ and $F\in C^{\infty
}(\mathbb{R}^{n})$, then $F(f_{1},...,f_{n})\in C^{\infty}(R).$

\item[2.3.] If $f:R\rightarrow\mathbb{R}$ is such that, for every $x\in R$,
there exist an open neighbourhood $U_{x}$ of $x$ and a function $f_{x}\in
C^{\infty}(R)$ satisfying
\[
f_{x}\mid U_{x}=f\mid U_{x},
\]
then $f\in C^{\infty}(R).$ Here the vertical bar $\mid$ denotes the restriction.
\end{description}

\noindent A differential space is a topological space endowed with a
differential structure. Clearly, smooth manifolds are differential spaces.

\begin{description}
\item [Lemma 1.]For every open subset $U$ of a differential space $R$ and
every $x\in U$, there exists $f\in C^{\infty}(R)$ satisfying $f\mid V=1$ for
some neighbourhood $V$ of $x$ contained in $U$, and $f\mid W=0$ for some open
subset $W$ of $R$ such that $U\cup W=R.$
\end{description}

\noindent\textbf{Proof} follows ref. \cite{sikorski}. Let $U$ be open in $R$
and $x\in U.$ It follows from condition 2.1 that there exists a map
$\varphi=(f_{1},...,f_{n}):R\rightarrow\mathbb{R}^{n}$, with $f_{1}%
,...,f_{n}\in C^{\infty}(R),$ and an open set $\tilde{U}\subseteq
\mathbb{R}^{n}$ such that $x\in\varphi^{-1}(\tilde{U})\subseteq U.$ Since
$\varphi(x)\in\tilde{U}\subseteq\mathbb{R}^{n}$, there exists $F\in C^{\infty
}(\mathbb{R}^{n})$ such that $F\mid\tilde{V}=1$ for some neigbourhood
$\tilde{V}$ of $\varphi(x)$ in $\mathbb{R}^{n}$ contained in $\tilde{U}$, and
$F\mid\tilde{W}=0$ for some open set $\tilde{W}$ in $R^{n}$ such that
$\tilde{U}\cup\tilde{W}=\mathbb{R}^{n}.$ Since $\varphi$ is continuous,
$V=\varphi^{-1}(\tilde{V})$ and $W=\varphi^{-1}(\tilde{W})$ are open in $V$.
Moreover, $\varphi^{-1}(\tilde{U})\subseteq U$ and $\tilde{U}\cup\tilde
{W}=\mathbb{R}^{n}$ imply that $U\cup W=R$. By condition 2.2, $f=F(f_{1}%
,...,f_{n})\in C^{\infty}(R)$. Furthermore, $f\mid V=F\raisebox{2pt}%
{$\scriptstyle\circ\,
$}\varphi\mid V=F\mid\varphi(V)=F\mid\tilde{V}=1.$ Similarly, $f\mid
W=F\mid\tilde{W}=0,$ which completes the proof. \hfill$\Box$

A continuous map $\varphi:S\rightarrow R$ between differential spaces $S$ and
$R$ is smooth if $\varphi^{\ast}f=f\raisebox{2pt}{$\scriptstyle\circ
\, $}\varphi\in C^{\infty}(S)$ for every $f\in C^{\infty}(R)$. A homeomorphism
$\varphi:S\rightarrow R$ is called a diffeomorphism if $\varphi$ and
$\varphi^{-1}$ are smooth.

If $R$ is a differential space with differential structure $C^{\infty}(R)$ and
$S$ is a subset of $R$, then we can define a differential structure
$C^{\infty}(S)$ on $S$ as follows. A function $f:S\rightarrow\mathbb{R}$ is in
$C^{\infty}(S)$ if and only if, for every $x\in S,$ there is an open
neighborhood $U$ of $x$ in $R$ and a function $f_{x}\in C^{\infty}(R)$ such
that $f|(S\cap U)=f_{x}|(S\cap U)$. The differential structure $C^{\infty}(S)$
described above is the smallest differential structure on $S$ such that the
inclusion map $\iota:S\rightarrow R$ is smooth. We shall refer to $S$ with the
differential structure $C^{\infty}(S)$ described above as a differential
subspace of $R$. If $S$ is a closed subset of $R,$ then the differential
structure $C^{\infty}(S)$ described above consists of restrictions to $S$ of
functions in $C^{\infty}(R)$.

A differential space $R$ is said to be locally diffeomorphic to a differential
space $S$ if, for every $x\in R$, there exists a neighbourhood $U$ of $x$
diffeomorphic to an open subset $V$ of $S$. More precisely, we require that
the differential subspace $U$ of $R$ be diffeomorphic to the differential
subspace $V$ of $S$. A differential space $R$ is a smooth manifold of
dimension $n$ if and only if it is locally diffeomorphic to $\mathbb{R}^{n}$.

Let $R$ be a differential space with a differential structure $C^{\infty}(R)
$. A derivation on $C^{\infty}(R)$ is a linear map $X:C^{\infty}(R)\rightarrow
C^{\infty}(R):f\mapsto X\cdot f$ satisfying Leibniz' rule
\begin{equation}
X\cdot(f_{1}f_{2})=(X\cdot f_{1})f_{2}+f_{1}(X\cdot f_{2}). \label{Leibniz}%
\end{equation}
We denote the space of derivations of $C^{\infty}(R)$ by $\mathrm{Der}%
C^{\infty}(R).$ It has the structure of a Lie algebra with the Lie bracket
$[X_{1},X_{2}]$ defined by
\[
\lbrack X_{1},X_{2}]\cdot f=X_{1}\cdot(X_{2}\cdot f)-X_{2}\cdot(X_{1}\cdot f)
\]
for every $X_{1},X_{2}\in\mathrm{Der}C^{\infty}(R)$ and $f\in C^{\infty}(R).$

\begin{description}
\item [Lemma 2.]If $f\in C^{\infty}(R)$ is a constant function, then $X\cdot
f=0$ for all $X\in\mathrm{Der}C^{\infty}(R).$
\end{description}

\noindent\textbf{Proof.} If $f\in C^{\infty}(R)$ is identically zero, then
$f^{2}=f=0$, and Leibniz' rule implies that $X\cdot f=X\cdot f^{2}=2f(X\cdot
f)=0$ for every $X\in\mathrm{Der}C^{\infty}(R)$. Similarly, if $f$ is a
non-zero constant function, that is $f(x)=c\neq0$ for all $x\in R$, then
$f^{2}=cf$, and the linearity of derivations implies that $X\cdot f^{2}%
=X\cdot(cf)=c(X\cdot f)$. On the other hand, Leibniz' rule implies that
$X\cdot f^{2}=2f(X\cdot f)=2c(X\cdot f)$. Hence $c(X\cdot f)=2c(X\cdot f)$.
Since $c\neq0$, it follows that $X\cdot f=0.$ \hfill$\Box$

\begin{description}
\item [Lemma 3.]If $f\in C^{\infty}(R)$ vanishes identically in an open set
$U\subseteq R,$ then $(X\cdot f)\mid U=0$ for all $X\in\mathrm{Der}C^{\infty}(R).$
\end{description}

\noindent\textbf{Proof}. If $f\in C^{\infty}(R)$ vanishes identically in an
open set $U\subseteq R$, then for each $x\in U$, there exists by Lemma 1 a
function $h\in C^{\infty}(R)$ such that $h(x)=1$ and $hf=0$. Therefore,
$0=X\cdot(hf)=h(X\cdot f)+f(X\cdot h)$ for every smooth derivation $X.$
Evaluating this identity at $x$, we get $(X\cdot f)(x)=0$ because $f(x)=0$.
Hence, $(X\cdot f)\mid U=0.$ \hfill$\Box$

\begin{description}
\item [Lemma 4.]Let $U$ be open in $R$, and $X_{U}$ a smooth derivation of
$C^{\infty}(U).$ For each $x\in U$, there exists an open neighbourhood $V$ of
$x$ contained in $U$, and $X\in\mathrm{Der}C^{\infty}(R)$ such that
\[
(X\cdot h)\mid V=(X_{U}\cdot(h\mid U))\mid V\text{ \quad for all \quad}h\in
C^{\infty}(R).
\]
\end{description}

\noindent\textbf{Proof}. Let $U$ be an open neighbourhood of $x_{0}$ in $R$,
and $X_{U}$ a smooth derivation of $C^{\infty}(U).$ There exist open sets $V$
and $W$ in $R$ such that $x_{0}\in V\subseteq\bar{V}\subseteq W\subseteq
\bar{W}\subseteq U$. Let $f\in C^{\infty}(R)$ be such that $f\mid\bar{V}=1$
and $f\mid S\backslash W=0$. Then $(f\mid U)X_{U}$ is a derivation of
$C^{\infty}(U)$ which vanishes on $U\backslash W$. Hence, it extends to a
smooth derivation $X$ of $C^{\infty}(R)$ such that, for every $h\in C^{\infty
}(R),$ $(X\cdot h)\mid V=(X_{U}\cdot(h\mid U))\mid V$.\hfill$\Box$

\bigskip

A local diffeomorphism $\varphi$ of $R$ to itself is a diffeomorphism
$\varphi:U\rightarrow V$, where $U$ and $V$ are open differential subspaces of
$R.$ For each $f\in C^{\infty}(R),$ the restriction of $f$ to $V$ is in
$C^{\infty}(V)$, and $\varphi^{\ast}f=f\raisebox{2pt}{$\scriptstyle\circ
\, $}\varphi$ is in $C^{\infty}(U)$. If $\varphi^{\ast}f$ coincides with the
restriction of $f$ to $U,$ we say that $f$ is $\varphi$-invariant, and write
$\varphi^{\ast}f=f.$ For each $X\in\mathrm{Der}(C^{\infty}(R)),$ the
restriction of $X$ to $U$ is in $\mathrm{Der}(C^{\infty}(U))$, and the
push-forward $\varphi_{\ast}X$ of $X$ by $\varphi$ is a derivation of
$C^{\infty}(V)$ such that
\begin{equation}
(\varphi_{\ast}X)\cdot(f\mid V)=\varphi^{-1\ast}(X\cdot(\varphi^{\ast
}f))\text{ for all }f\in C^{\infty}(R). \label{push-forward}%
\end{equation}
Since all functions in $C^{\infty}(V)$ locally coincide with restrictions to
$V$ of functions in $C^{\infty}(R)$, equation (\ref{push-forward}) determines
$\varphi_{\ast}X$ uniquely. If $\varphi_{\ast}X$ coincides with the
restriction of $X$ to $V$, we say that $X$ is $\varphi$-invariant and write
$\varphi_{\ast}X=X.$

\section{Subcartesian spaces.}

Subcartesian spaces were introduced by Aronszajn, \cite{aronszajn}, and
developed in \cite{aronszajn-szeptycki} and \cite{marshall}. They are
Hausdorff differential spaces locally diffeomorphic to a differential subspace
of a Cartesian space. In other words, a subcartesian space is a Hausdorff
differential space $S$ that can be covered by open sets, each of which is
diffeomorphic to a differential subspace of a Cartesian space. In the
remainder of this paper, we restrict our considerations to differential spaces
that are Hausdorff, second countable and paracompact.

In this section, we describe properties of differential subspaces of
$\mathbb{R}^{n}$ which extend to subcartesian spaces. In the remainder of this
section, $R$ denotes a differential subspace of $\mathbb{R}^{n},$ considered
as a differential space endowed with the standard differential structure
$C^{\infty}(\mathbb{R}^{n}).$ In other words, a function $f:R\rightarrow
\mathbb{R} $ is in $C^{\infty}(R)$ if, for every $x\in R$, there exists an
open set $U$ in $\mathbb{R}^{n}$ and $f_{U}\in C^{\infty}(\mathbb{R}^{n})$
such that $f\mid U\cap R=f_{U}\mid U\cap R.$

\begin{description}
\item [Lemma 5.]Let $W$ be an open subset of $R\subseteq\mathbb{R}^{n}$, and
$f_{W}\in C^{\infty}(W)$. For every $x\in W$ there exists a function $f\in
C^{\infty}(R)$ and a neighbourhood $V$ of $x$ contained in $W$ such that
$f\mid V=f_{W}\mid V.$
\end{description}

\noindent\textbf{Proof. }The proof is an immediate consequence of the
definition of a differential subspace. \hfill$\Box$

For every differential space $S,$ each $X\in\mathrm{Der}C^{\infty}(S)$ and
every $x\in S$, we denote by $X(x):C^{\infty}(S)\rightarrow\mathbb{R}$ the
composition of the derivation $X$ with the evaluation at $x$. In other words,
$X(x)\cdot f=(X\cdot f)(x)$ for all $f\in C^{\infty}(S)$. We use the notation
\begin{equation}
\mathrm{Der}_{x}C^{\infty}(S)=\{X(x)\mid X\in\mathrm{Der}C^{\infty}(S)\},
\label{Der_x}%
\end{equation}
and refer to elements of $\mathrm{Der}_{x}C^{\infty}(S)$ as derivations of
$C^{\infty}(S)$ at $x$.

\begin{description}
\item [Lemma 6.]Let $R$ be a differential subspace of $\mathbb{R}^{n}$ and $X$
a derivation of $C^{\infty}(R).$ For every $x\in R,$ there is $\tilde{X}%
(x)\in\mathrm{Der}_{x}C^{\infty}(\mathbb{R}^{n})$ such that
\begin{equation}
\tilde{X}(x)\cdot f=X(x)\cdot(f\mid R)\text{ \quad for all \quad\ }f\in
C^{\infty}(\mathbb{R}^{n}). \label{derivation x}%
\end{equation}
\end{description}

\noindent\textbf{Proof}. Let $\tilde{X}(x):C^{\infty}(\mathbb{R}%
^{n})\rightarrow\mathbb{R}$ be given by equation (\ref{derivation x}). It is a
linear map satisfying Leibniz' rule $\tilde{X}(x)\cdot(fh)=f(x)\left(
\tilde{X}(x)\cdot h\right)  +h(x)\left(  \tilde{X}(x)\cdot f\right)  $. Hence,
$\tilde{X}(x)\in T_{x}\mathbb{R}^{n},$ and it extends to a smooth vector field
$\tilde{X}$ on $\mathbb{R}^{n}$ so that $\tilde{X}(x)$ is the value of
$\tilde{X} $ at $x\in\mathbb{R}^{n}$. Since vector fields on $\mathbb{R}^{n} $
are derivations of $C^{\infty}(\mathbb{R}^{n})$, it follows that $\tilde
{X}(x)\in\mathrm{Der}_{x}C^{\infty}(\mathbb{R}^{n}).$ \hfill$\Box$

\begin{description}
\item [Proposition 1.]Let $R$ be a differential subspace of $\mathbb{R}^{n}$
and $X$ a derivation of $C^{\infty}(R).$ For every $f_{1},...,f_{m}\in
C^{\infty}(R)$ and every $F\in C^{\infty}(\mathbb{R}^{m}),$%
\begin{equation}
X\cdot F(f_{1},...,f_{m})=\sum_{i=1}^{m}\partial_{i}F(f_{1},...,f_{m})(X\cdot
f_{i}). \label{smooth derivation}%
\end{equation}
\end{description}

\noindent\textbf{Proof. }Let $f_{1},...,f_{m}\in C^{\infty}(R)$ and $x\in R $.
We denote by\textbf{\ }$p_{1},...,p_{n}:\mathbb{R}^{n}\rightarrow\mathbb{R}$
the coordinate functions on $\mathbb{R}^{n}$. There exists a neighbourhood $U$
of $x\in\mathbb{R}^{n}$ and functions $F_{1},...,F_{m}\in C^{\infty
}(\mathbb{R}^{n})$ such that $f_{i}\mid U\cap R=F_{i}(p_{1},...,p_{n})\mid
U\cap R.$ Hence, for every $F\in C^{\infty}(\mathbb{R}^{m})$,
\[
F(f_{1},...,f_{m})\mid U\cap R=F(F_{1}(p_{1}...,p_{n}),...,F_{m}%
(p_{1},...,p_{n}))\mid U\cap R\text{.}%
\]
By Lemma 4, there exists $\tilde{X}(x)\in\mathrm{Der}_{x}C^{\infty}%
(\mathbb{R}^{n})$ such that equation (\ref{derivation x}) is satisfied.
Hence,
\begin{align*}
&  X(x)\cdot F(f_{1},...,f_{m})=X(x)\cdot\left(  F(F_{1}(p_{1}...,p_{n}%
),...,F_{m}(p_{1},...,p_{n}))\mid R\right) \\
&  =\tilde{X}(x)\cdot\left(  F(F_{1}(p_{1}...,p_{n}),...,F_{m}(p_{1}%
,...,p_{n}))\right) \\
&  =\sum_{i=1}^{n}\left(  \partial_{i}F(F_{1}(p_{1}...,p_{n}),...,F_{m}%
(p_{1},...,p_{n}))\right)  (x)\left(  \tilde{X}(x)\cdot F_{i}(p_{1}%
...,p_{n})\right) \\
&  =\sum_{i=1}^{n}\left(  \partial_{i}F(F_{1}(p_{1}...,p_{n}),...,F_{m}%
(p_{1},...,p_{n}))\right)  (x)\left(  X(x)\cdot\left(  F_{i}(p_{1}%
...,p_{n})\mid R\right)  \right) \\
&  =\sum_{i=1}^{n}\left(  \partial_{i}F(f_{1},...,f_{m})\right)  (x)\left(
X(x)\cdot f_{i}\right) \\
&  =\sum_{i=1}^{n}\left(  \partial_{i}F(f_{1},...,f_{m})\right)  (x)\left(
X\cdot f_{i}\right)  (x).
\end{align*}
This holds for every $x\in R$, which implies equation (\ref{smooth
derivation}). \hfill$\Box$

\begin{description}
\item [Proposition 2.]Let $R$ be a differential subspace of $\mathbb{R}^{n}. $
For every $X\in\mathrm{Der}C^{\infty}(R)$ and each $x\in R$, there exists an
open neighbourhood $U$ of $x$ in $\mathbb{R}^{n}$ and $\tilde{X}%
\in\mathrm{Der}C^{\infty}(U)$ such that
\[
\left(  \tilde{X}\cdot\left(  f\mid U\right)  \right)  \mid U\cap R=\left(
X\cdot\left(  f\mid R\right)  \right)  \mid U\cap R\text{ \quad for
all\quad\ }f\in C^{\infty}(\mathbb{R}^{n}).
\]
\end{description}

\noindent\textbf{Proof}. Let $h_{1},...,h_{n}$ be the restrictions to $R$ of
Cartesian coordinates $p_{1},...,p_{n}$ on $\mathbb{R}^{n}$. For every
$X\in\mathrm{Der}(C^{\infty}(R)),$ the functions $X\cdot h_{1},...,X\cdot
h_{n}$ are in $C^{\infty}(R).$ Hence, for every $x\in R$, there exists an open
neighbourhood $U$ of $x$ in $\mathbb{R}^{n}$ and functions $f_{1},...,f_{n}\in
C^{\infty}(\mathbb{R}^{n})$ such that $(X\cdot h_{i})\mid U\cap R=f_{i}\mid
U\cap R$ for $i=1,...,n.$

Let $f\in C^{\infty}(\mathbb{R}^{n})$. Then $f\mid U\cap R=f(h_{1}%
,...,h_{n})\mid U\cap R$ and
\begin{align*}
\left(  X\cdot(f\mid R)\right)  \mid U\cap R  &  =\left(  X\cdot
(f(h_{1},...,h_{n}))\right)  \mid U\cap R\\
&  =\sum_{i=1}^{n}\left(  \partial_{i}f(h_{1},...,h_{m})\mid U\cap R\right)
(X\cdot h_{i})\mid U\cap R\\
&  =\sum_{i=1}^{n}\left(  (\partial_{i}f)\mid U\cap R\right)  \left(
f_{i}\mid U\cap R\right)  .
\end{align*}
Consider the vector field
\begin{equation}
\tilde{X}=f_{1}\partial_{1}+...+f_{n}\partial_{n} \label{Xtilde}%
\end{equation}
on $U$. Since $U$ is open in $\mathbb{R}^{n}$, $\tilde{X}$ is a derivation of
$C^{\infty}(U)$. Moreover,
\[
\left(  X\cdot(f\mid R)\right)  \mid U\cap R=\left(  \tilde{X}\cdot\left(
f\mid U\right)  \right)  \mid U\cap R
\]
for all $f\in C^{\infty}(\mathbb{R}^{n}).$ \hfill$\Box$

\bigskip

It follows from Proposition 2 that every derivation $X$ of $C^{\infty}(R)$ can
be locally extended to a derivation of $C^{\infty}(\mathbb{R}^{n})$. Clearly,
this extension need not be unique. Moreover, If $\tilde{X}$ is a smooth vector
field on $\mathbb{R}^{n}$, and $f\in C^{\infty}(\mathbb{R}^{n}),$ then the
restriction of $\tilde{X}\cdot f$ to $R$ need not be determined by the
restriction of $f$ to $R$. Hence, not every vector field on $\mathbb{R}^{n}$
restricts to a derivation of $C^{\infty}(R).$

\begin{description}
\item [Lemma 7.]Let $U$ be an open subset of $R\subseteq\mathbb{R}^{n},$ and
$X_{U}$ a smooth derivation of $C^{\infty}(U)$. For each $x\in U$, there
exists an open neighbourhood $V$ of $x$ contained in $U$, and $X\in
\mathrm{Der}C^{\infty}(R)$ such that
\[
(X\cdot h)\mid V=(X_{U}\cdot(h\mid U))\mid V\text{ \quad for all}\quad\text{
}h\in C^{\infty}(R).
\]
\end{description}

\noindent\textbf{Proof}. By Lemma 1, there exist open sets $V$ and $W$ in $R,$
such that $x\in V\subseteq U$ and $U\cup W=R,$ and $f\in C^{\infty}(R)$
satisfying $f\mid V=1$ and $f\mid W=0.$ For every $h\in C^{\infty}(R)$, let
$X$ be given by
\begin{align*}
(X\cdot h)  &  \mid U=(f\mid U)X_{U}\cdot(h\mid U)\\
(X\cdot h)  &  \mid W=0.
\end{align*}
In $U\cap W$ we have
\[
\left(  (f\mid U)X_{U}\cdot(h\mid U)\right)  \mid U\cap W=(f\mid U\cap
W)\left(  X_{U}\cdot(h\mid U)\right)  \mid U\cap W=0
\]
because $f\mid W=0$. Hence, $X$ is well defined. Moreover, $X\in
\mathrm{Der}C^{\infty}(R),$ and $f\mid V=1$ implies $(X\cdot h)\mid
V=(X_{U}\cdot(h\mid U))\mid V$ for all $h\in C^{\infty}(R)$. \hfill$\Box$

\begin{description}
\item [Lemma 8.]Let $R$ be a differential subspace of $\mathbb{R}^{n}$. If $U$
and $V$ are open subsets of $\mathbb{R}^{n}$ and $\varphi:U\rightarrow V$ is a
diffeomorphism such that $\varphi(U\cap R)=V\cap R$, then the restriction of
$\varphi$ to $U\cap R$ is a diffeomorphism of $U\cap R$ onto $V\cap R.$
\end{description}

\noindent\textbf{Proof}. By assumption, $R$ is a topological subspace of
$\mathbb{R}^{n}$, the mapping $\varphi:U\rightarrow V$ is a homeomorphism, and
$\varphi(U\cap R)=V\cap R$. Hence, for every open subset $W$ of $\mathbb{R}%
^{n}, $ $\varphi^{-1}(W\cap(V\cap R))$ is open in $U\cap R$, and
$\varphi(W\cap(U\cap R))$ is open in $V\cap R$. Thus, $\varphi$ induces a
homeomorphism $\psi:U\cap R\rightarrow V\cap R$.

Moreover, $\varphi$ induces a diffeomorphism of open differential subspaces
$U$ and $V$ of $\mathbb{R}^{n}.$ We want to show that $f\in C^{\infty}(V\cap
R) $ implies $\psi^{\ast}f\in C^{\infty}(U\cap R).$ Given $x\in U\cap R$, let
$y=\psi(x)\in V\cap R$. Since $R$ is a differential subspace of $\mathbb{R}%
^{n}$ and $f\in C^{\infty}(V\cap R)$, there exists a neighbourhood $W$ of
$x\in V$ and a function $f_{W}\in C^{\infty}(V)$ such that $f\mid W\cap
R=f_{W}\mid W\cap R$. Moreover, $\varphi^{-1}(W)$ is a neighbourhood of $x$ in
$U$, $\varphi^{\ast}f_{W}$ is in $C^{\infty}(U)$, and
\begin{align*}
\psi^{\ast}f\mid(\varphi^{-1}(W)\cap R)  &  =f\raisebox{2pt}{$\scriptstyle
\circ\, $}\psi\mid(\varphi^{-1}(W)\cap R)=f\raisebox{2pt}{$\scriptstyle
\circ\, $}(\varphi\mid(\varphi^{-1}(W)\cap R))\\
&  =f\mid(W\cap R)=f_{W}\mid(W\cap R)\\
&  =f_{W}\raisebox{2pt}{$\scriptstyle\circ\, $}(\varphi\mid(\varphi
^{-1}(W)\cap R))=f_{W}\raisebox{2pt}{$\scriptstyle\circ\, $}\varphi
\mid(\varphi^{-1}(W)\cap R)\\
&  =\varphi^{\ast}f_{W}\mid(\varphi^{-1}(W)\cap R).
\end{align*}
Thus, for every $x\in U\cap R$, there exists a neighbourhood $\varphi^{-1}(W)
$ of $x$ in $U$ and a function $\varphi^{\ast}f_{W}$ in $C^{\infty}(U)$ such
that $\psi^{\ast}f\mid(\varphi^{-1}(W)\cap R)=\varphi^{\ast}f_{W}\mid
(\varphi^{-1}(W)\cap R).$ This implies that $\psi^{\ast}f\in C^{\infty}(U\cap R).$

It follows that $\psi$ is smooth. In a similar manner we can prove that
$\psi^{-1}$ is smooth. Hence, $\psi$ is a diffeomorphism. \hfill$\Box\bigskip$

Let $S$ be a subcartesian space. It can be covered by open sets, each of which
is diffeomorphic to a differential subspace $R$ of $\mathbb{R}^{n}$. All the
properties of differential subspaces of $\mathbb{R}^{n}$ discussed in Lemmata
4 through 8 and Propositions 1 and 2 are local. Hence they extend to
subcartesian spaces.

\section{Families of Vector Fields}

In this section, we discuss properties of vector fields on subcartesian spaces.

In the category of smooth manifolds, translations along integral curves of a
smooth vector field give rise to local diffeomorphisms. In the category of
differential spaces, not all derivations generate local diffeomorphisms. We
reserve the term vector field for a derivation that generates a local
one-parameter group of local diffeomorphisms.

Let $I$ be an interval in $\mathbb{R}$. A smooth curve $c:I\rightarrow S$ on a
differential space $S$ is an integral curve of\textit{\ }$X\in\mathrm{Der}%
(C^{\infty}(S))$ if
\[
(X\cdot f)(c(t))=\frac{d}{dt}f(c(t))
\]
for every $f\in C^{\infty}(S)$ and $t\in I$. If $0\in I$ and $c(0)=x,$ we
refer to $c$ as an integral curve of $X$\thinspace through $x.$

\begin{description}
\item [Theorem 1.]Assume that $S$ is a subcartesian space. For every $x\in S $
and every $X\in\mathrm{Der}(C^{\infty}(S))$ there exists a unique maximal
integral curve $c:I\rightarrow S$ through $x$.
\end{description}

\noindent\textbf{Proof Outline} (a detailed proof is given in \cite{sniatycki
2002}). Since $S$ is a subcartesian space, given $x\in S$, there exists a
neighbourhood $V$ of $x$ in $S$ diffeomorphic to a differential subspace $R$
of $\mathbb{R}^{n}.$ In order to simplify the notation, we use the
diffeomorphism between $R$ and $V$ to identify them, and write $V=R$. By
Proposition 2, there exists an extension of $X$ to a smooth vector field
$\tilde{X}$ on an open neighbourhood $U$ of $x$ in $\mathbb{R}^{n}$ given by
equation (\ref{Xtilde}).

Given $y\in R\subseteq\mathbb{R}^{n}$, consider an integral curve $\tilde
{c}:\tilde{I}\rightarrow\mathbb{R}^{n}$ of $\tilde{X}$ such that $\tilde
{c}(0)=y$. Let $I$ be the connected component of $\tilde{c}^{-1}(R)$
containing $0,$ and $c:I\rightarrow R$ the curve in $R$ obtained by the
restriction of $\tilde{c}$ to $I$. Then, $c(0)=y$. Moreover, for each $t\in
I$, and $f\in C^{\infty}(S)$ there exists a neigbourhood $U$ of $\tilde{c}(t)$
in $\mathbb{R}^{n}$ and a function $F\in C^{\infty}(\mathbb{R}^{n})$ such that
$f\mid R\cap U=F\mid R\cap U$. Hence, by Proposition 2,
\[
\frac{d}{dt}f(c(t))=\frac{d}{dt}F(\tilde{c}(t))=(\tilde{X}\cdot F)(\tilde
{c}(t))=(X\cdot f)(c(t)).
\]
This implies that $c:I\rightarrow R$ is an integral curve of $X$ through $y.$
Since $I$ is a connected subset of $\mathbb{R}$, it is an interval. Local
uniqueness of $c$ (up to an extension of the domain) follows from the local
uniqueness of solutions of differential equations on $\mathbb{R}^{n}.$

The above argument gives existence and local uniqueness of integral curves of
derivations of $C^{\infty}(S).$ The usual technique of patching local
solutions, and the fact that the union of intervals with pairwise non-empty
intersection is an interval, lead to the global uniqueness of integral curves
of derivations on a subcartesian space $S$.\hfill$\Box\bigskip$

Let $X$ be a derivation of $C^{\infty}(S)$. We denote by $\varphi_{t}(x)$, the
point on the maximal integral curve of $X$ through $x$ corresponding to the
value $t$ of the parameter. Given $x\in\mathbb{R}^{n}$, $\varphi_{t}(x) $ is
defined for $t$ in an interval $I_{x}$ containing zero, and $\varphi_{0}%
(x)=x$. If $t,$ $s$ and $t+s$ are in $I_{x},$ $s\in I_{\varphi_{t}(x)}$ and
$t\in I_{\varphi_{s}(x)},$ then $\varphi_{t+s}(x)=\varphi_{s}(\varphi
_{t}(x))=\varphi_{t}(\varphi_{s}(x))$. In the case when $S$ is a manifold, the
map $\varphi_{t}$ is a diffeomorphism of a neighbourhood of $x$ in $S$ onto a
neighbourhood of $\varphi_{t}(x)$ in $S$. For a subcartesian space $S,$ the
map $\varphi_{t}$ might fail to be a local diffeomorphism. We adopt the following

\begin{description}
\item [Definition of a vector field.]Let $S$ be a subcartesian space. A
derivation $X$ of $C^{\infty}(S)$ is a vector field on $S$ if, for every $x\in
S$, there exists an open neighbourhood $U$ of $x$ in $S$, and $\varepsilon>0$
such that, for every $t\in(-\varepsilon,\varepsilon),$ the map $\varphi_{t}$
is defined on $U$, and its restriction to $U$ is a diffeomorphism from $U$
onto an open subset of $S$.
\end{description}

\noindent\textbf{Example 1.} Consider $S=[0,\infty)\subseteq\mathbb{R}$ with
the structure of a differential subspace of $\mathbb{R}.$ Let $(X\cdot
f)=f^{\prime}(x)$ for every $f\in C^{\infty}([0,\infty))$ and $x\in
\lbrack0,\infty).$ Note that the derivative at $x=0$ is is the right
derivative; it is uniquely defined by $f(x)$ for $x\geq0$. For this $X$, the
map $\varphi_{t}$ is given by $\varphi_{t}(x)=x+t$ whenever $x$ and $x+t$ are
in $[0,\infty)$. In particular, for every neighbourhood $U$ of $0$ in
$[0,\infty)$ there exists $\delta>0$ such that $[0,\delta)\subseteq U$.
Moreover, $\varphi_{t}$ maps $[0,\delta)$ onto $[t,\delta+t)$ which is not a
neighbourhood of $t=\varphi_{t}(0)$ in $[0,\infty)$. Hence, the derivation $X$
is not a vector field on $[0,\infty)$. On the other hand, for every $f\in
C^{\infty}[0,\infty)$ such that $f(0)=0$, the derivation $fX$ is a vector
field, because $0$ is a fixed point of $fX$.\bigskip

We say that a subcartesian space $S$ is locally closed if every point of $S$
has a neighbourhood diffeomorphic to the intersection of an \ open and a
closed subset of a Cartesian space. There is a simple criterion characterizing
vector fields on a locally closed subcartesian space; namely,

\begin{description}
\item [Proposition 3.]Let $S$ be a locally closed subcartesian space. A
derivation $X$ of $C^{\infty}(S)$ is a vector field on $S$ if the domain of
every maximal integral curve of $X$ is open in $\mathbb{R}$.
\end{description}

\noindent\textbf{Proof}. Consider first the case when $S$ is a relatively
closed differential subspace of $\mathbb{R}^{n}.$ In other words, $S=U\cap C$,
where $U$ is open and $C$ is closed in $\mathbb{R}^{n}.$ By Proposition 2, we
may assume that $X$ on $S$ extends to a vector field $\tilde{X}$ on $U$. We
denote by $\tilde{\varphi}_{t}$ the local one-parametr group of local
diffeomorphisms of $U$ generated by $\tilde{X}$. By Theorem 1, for every $x\in
S$, there is a maximal interval $I_{x}\in\mathbb{R}$ such that $\varphi
_{t}(x)=\tilde{\varphi}_{t}(x)\in S$ for all $t\in I_{x}$.

Given $x\in S\subseteq U$, there exists $\varepsilon>0$ and a neighbourhood
$U^{\prime}$ of $x$ in $U$ such that such that, for every $t\in(-\varepsilon
,\varepsilon),$ the map $\tilde{\varphi}_{t}$ is defined on $U^{\prime}$, and
its restriction to $U^{\prime}$ is a diffeomorphism from $U^{\prime}$ onto an
open subset of $U.$ In view of Lemma 8, it suffices to show that there exists
$\delta\in(0,\varepsilon]$ and a neighbourhood $U^{\prime\prime}$ of $x$ in
$U^{\prime}$ such that $\tilde{\varphi}_{t}$ maps $U^{\prime\prime}\cap C$ to
$\tilde{\varphi}_{t}(U^{\prime\prime})\cap C$ for all $t\in(-\delta,\delta).$
Suppose that there are no $U^{\prime\prime}$ and $\delta$ satisfying this
condition. This means that, for every neighbourhood $U^{\prime\prime}$ of $x$
in $U$, and every $\delta\in(0,\varepsilon]$, there exists a point $y\in
U^{\prime\prime}\cap C$ and $s\in(-\delta,\delta)$ such that $\tilde{\varphi
}_{s}(y)\notin\tilde{\varphi}_{s}(U^{\prime\prime})\cap C$. Since
$\tilde{\varphi}_{t}(y)\in\tilde{\varphi}_{t}(U^{\prime\prime})$ for every
$t\in(-\varepsilon,\varepsilon)$, it follows that $\tilde{\varphi}%
_{s}(y)\notin C$. Hence, $s$ is not in the domain $I_{y}$ the maximal integral
curve of $X$ through $y$. If $s>0$, let $u$ be the infimum of the set
$\{t\in\lbrack0,s]\mid\tilde{\varphi}_{t}(y)\notin C\}.$ Then, $\varphi
_{t}(y)\in C$ for all $t\in\lbrack0,u)$. Since $\varphi_{t}(y)$ is continuous
in $t$ and $C$ is closed, it follows that $\varphi_{u}(y)\in C$. Moreover, for
every $v>u$, there exists $t\in(u,v)$ such that $\varphi_{t}(y)\notin C$. It
implies that $[0,\infty)\cap I_{y}=[0,u]$. Hence, the domain $I_{y}$ of the
maximal integral curve of $X$ through $y$ is not open in $\mathbb{R}$,
contrary to the assumption of the theorem. Hence, the case $s>0$ is excluded.
Similarly, we can show that the case $s<0$ is inconsistent with the assumption
that the domains of all maximal integral curves of $X$ are open.

We have shown that there exists $\delta\in(0,\varepsilon]$ and a neighbourhood
$U^{\prime\prime}$ of $x$ in $U^{\prime}$ such that $\tilde{\varphi}_{t}$ maps
$U^{\prime\prime}\cap C$ to $\tilde{\varphi}_{t}(U^{\prime\prime})\cap C$ for
all $t\in(-\delta,\delta).$ This implies that $\varphi_{t}(z)=\tilde{\varphi
}_{t}(z)$ is defined for every $t\in(-\delta,\delta),$ and each $z\in
U^{\prime\prime}$. By Lemma 8, it follows that, $\varphi_{t}$ restricted to
$U^{\prime\prime}\cap S$ is a diffeomorphism onto $\varphi_{t}(U^{\prime
\prime})\cap S$. Since this holds for every $x\in S$, we conclude that
$\varphi_{t}$ is a local one-parameter group of local diffeomorphisms of $S$.
Hence, $X$ is a vector field on $S$.

Consider now a derivation $X$ on a locally closed subcartesian space $S$ such
that the domains of all maximal integral curves of $X$ are open. For each
$x\in S$, we denote by $\varphi_{t}(x)$ the point on the maximal integral
curve of $X$ through $x$ corresponding to the value $t$ of the parameter. The
function $t\mapsto\varphi_{t}(x)$ is defined on an interval $I_{x}$ in
$\mathbb{R}$, which is open by hypothesis.

For every $x\in S$ there exists a neighbourhood $W$ of $x$ in $S$ and a
diffeomorphism $\chi$ of $W$ onto a locally closed subspace $U\cap C$ of
$\mathbb{R}^{n}$. By the first part of the proof, the push-forward of $X$ by
$\chi$ is a vector field on $U\cap C$. Since $\chi$ is a diffeomorphism, it
follows that there exists a neighbourhood $W^{\prime}$ of $x$ in $W\subseteq
S$ and $\varepsilon>0$ such that, for every $t\in(-\varepsilon,\varepsilon),$
the map $\varphi_{t}$ is defined on $W^{\prime}$, and its restriction to
$W^{\prime}$ is a diffeomorphism from $W^{\prime}$ onto an open subset of
$W\subseteq S$. Hence, $X$ is a vector field on $S$. \hfill$\Box$\bigskip

The following example shows that the assumption that $S$ is locally closed is
essential in Proposition 2.\bigskip

\noindent\textbf{Example 2}. The set
\[
S=\{(x_{1},x_{2})\in\mathbb{R}^{2}\mid x_{1}^{2}+(x_{2}-1)^{2}<1\text{ or
}x_{2}=0\}
\]
is not locally closed. The vector field $\tilde{X}=\frac{\partial}{\partial
x_{1}}$ on $\mathbb{R}^{2}$ restricts to a derivation $X$ of $C^{\infty}(S)$.
For every $x=(x_{1},x_{2})\in S$, $\tilde{\varphi}_{t}(x)=(x_{1}+t,x_{2})$ for
all $t\in\mathbb{R}$. Its restriction to $S$ induces $\varphi_{t}$ given by
$\varphi_{t}(x_{1},x_{2})=(x_{1}+t,x_{2})$ for $t\in\left(  -x_{1}%
-\sqrt{1-(x_{2}-1)^{2}},-x_{1}+\sqrt{1-(x_{2}-1)^{2}}\right)  $ if $x_{2}>0,$
and for $t\in\mathbb{R}$ if $x_{2}=0.$ Hence, all integral curves of $X$ have
open domains. Nevertheless, $\varphi_{t}$ fails to be a local one-parameter
local group of diffeomorphisms of $S$. \bigskip

Let $\mathcal{F}$ be a family of vector fields on a subcartesian space $S$.
For each $X\in\mathcal{F}$, we denote by $\varphi_{t}^{X}$ the local
one-parameter group of local diffeomorphisms of $S$ generated by $X$. We say
that the family $\mathcal{F}$ is locally complete if, for every $X,Y\in
\mathcal{F}$, $t\in\mathbb{R}$ and $x\in S,$ for which $\varphi_{t\ast}%
^{X}Y(x)$ is defined, there exists an open neighbourhood $U$ of $x$ and
$Z\in\mathcal{F}$ such that $\varphi_{t\ast}^{X}Y\mid U=Z\mid U.$

For example, a family consisting of a single vector field $X$ is locally
complete because $\varphi_{t\ast}^{X}X(x)=X(x)$ at all points $x\in S$ for
which $\varphi_{t}(x)$ is defined.

\begin{description}
\item [Theorem 2.]The family $\mathcal{X}(S)$ of all vector fields on a
subcartesian space $S$ is locally complete.
\end{description}

\noindent\textbf{Proof}. For $X\in\mathcal{X}(S)$, let $\varphi_{t}^{X}$
denote the local one-parameter group of local diffeomorphisms of $S$ generated
by $X$. For a given $t\in\mathbb{R}$ let $U$ be a neighbourhood of $x\in S$
such that $\varphi_{t}^{X}$ maps $U$ diffeomorphically onto an open subset $V$
of $S$. For each $Y\in\mathcal{X}(S)$, $\varphi_{t\ast}^{X}Y$ is in
$\mathrm{Der}(C^{\infty}(V)).$ If $\varphi_{s}^{Y}$ denotes the local
one-parameter group of diffeomorphisms of $Y$, then
\begin{align*}
\frac{d}{ds}f(\varphi_{t}^{X}(\varphi_{s}^{Y}(x))  &  =\frac{d}{ds}%
(\varphi_{t}^{X\ast}f)(\varphi_{s}^{Y}(x))=(Y\cdot(\varphi_{t}^{X\ast
}f))(\varphi_{s}^{Y}(x))\\
&  =(\varphi_{t}^{X\ast}\left(  \varphi_{t\ast}^{X}Y\cdot f\right)
)(\varphi_{s}^{Y}(x))=\left(  \varphi_{t\ast}^{X}Y\cdot f\right)  (\varphi
_{t}^{X}(\varphi_{s}^{Y}(x)))
\end{align*}
for every $f\in C^{\infty}(V)$ and $x\in V$ and $s\in I_{y}$ such that
$\varphi_{s}^{Y}(x)\in U$. Hence, $s\mapsto\varphi_{t}^{X}(\varphi_{s}%
^{Y}(x))$ is an integral curve of $\varphi_{t\ast}^{X}Y$ through $\varphi
_{t}^{X}(x)$. Since $Y$ is a vector field, for every $x\in S$, there exists an
open neighbourhood $W$ of $x$ in $S$, and $\varepsilon>0$ such that, for every
$s\in(-\varepsilon,\varepsilon),$ the map $\varphi_{s}^{Y}$ is defined on $W$,
and its restriction to $W$ is a diffeomorphism from $W$ onto an open subset of
$S$. We can choose $W$ and $\varepsilon$ so that $\varphi_{s}^{Y}$ maps $W$
into $U$ for all $s\in(-\varepsilon,\varepsilon)$. Since $\varphi_{t}^{X}%
\,\ $maps $W$ diffeomorphically ontu $V$, it follows that $\varphi_{t}%
^{X}\raisebox{2pt}{$\scriptstyle\circ\, $}\varphi_{s}^{Y}$ restricted to $W$
maps $W$ difeomorphically onto $\varphi_{t}^{X}(\varphi_{s}^{Y}(W)),$ which is
open in $V$ for all $s\in$ $(-\varepsilon,\varepsilon).$ Hence $\varphi
_{t\ast}^{X}Y$ is a vector field on $V$.

For every $x\in V$, there exists an open neighbourhood $W$ of $x$ such that
$\bar{W}\subset V$. Let $f\in C^{\infty}(V)$ be such that $f(x)=1$ and $f$
vanishes identically on $V\backslash W$. Then $f\varphi_{t*}Y$ is a vector
field on $V$ vanishing on $V\backslash W,$ and it extends to a vector field
$Z$ on $S$. Hence, $(\varphi_{t*}Y)(x)=f(x)(\varphi_{t*}Y)(x)=Z(x)$ and
$Z\in\mathcal{X}(S)$.

The above argument is valid for every $X$ and $Y$ in $\mathcal{X}(S)$. Hence,
$\mathcal{X}(S)$ is a locally complete family of vector fields. \hfill$\Box$

\section{Orbits and integral manifolds}

In this section we prove that orbits of families of vector fields on a
subcartesian space $S$ are manifolds. This is an extension of the results of
Sussmann, \cite{sussmann}, to the category of subcartesian spaces.

Let $\mathcal{F}$ be a family of vector fields on a subcartesian space $S$.
For each $X\in\mathcal{F}$ we denote by $\varphi_{t}^{X}$ the local
one-parameter group of diffeomorphisms of $S$ generated by $X$. The family
$\mathcal{F}$ gives rise to an equivalence relation $\sim$ on $S$ defined as
follows: $x\sim y$ if there exist vector fields $X^{1},...X^{n}\in\mathcal{F}$
and $t_{1},...,t_{n}\in\mathbb{R}$ such that
\begin{equation}
y=\left(  \varphi_{t_{n}}^{X^{n}}\raisebox{2pt}{$\scriptstyle\circ
\, $}...\raisebox{2pt}{$\scriptstyle\circ\, $}\varphi_{t_{1}}^{X^{1}}\right)
(x). \label{orbit}%
\end{equation}
In other words, $x\sim y$ if there exists a piecewise smooth curve $c$ in $S,
$ with tangent vectors given by restrictions to $c$ of vector fields in
$\mathcal{F},$ which joins $x$ and $y$. The equivalence class of this relation
containing $x$ is called the orbit of $\mathcal{F}$ through $x.$ The aim of
this section is to prove that orbits of locally complete families of vector
fields on a subcartesian space $S,$ are manifolds and give rise to a singular
foliation of $S$. This is an extension of the results of Sussmann,
\cite{sussmann}, to the category of differential spaces.

Following Sussmann's notation, we write $\xi=(X^{1},...,X^{m})$,
$T=(t_{1},...,t_{m})$ and
\[
\xi_{T}(x)=\left(  \varphi_{t_{m}}^{X^{m}}\raisebox{2pt}{$\scriptstyle\circ
\, $}...\raisebox{2pt}{$\scriptstyle\circ\, $}\varphi_{t_{1}}^{X^{1}}\right)
(x).
\]
The expression for $\xi_{T}(x)$ is defined for all $(T,x)$ in an open subset
$\Omega(\xi)$ of $\mathbb{R}^{m}\times S.$ Let $\Omega_{T}(\xi)$ denote the
set of all $x\in S$ such that $(T,x)\in\Omega(\xi)$. In other words,
$\Omega_{T}(\xi)$ is the set of all $x$ for which $\xi_{T}(x)$ is defined.
Moreover, we denote by $\Omega_{\xi,x}\subseteq\mathbb{R}^{m}$ the set of
$T\in\mathbb{R}^{m}$ such that $\xi_{T}(x)$ is defined.

We now assume that $S$ is a subset of $\mathbb{R}^{n}$. For each $x\in
S\subseteq\mathbb{R}^{n}$ and $\xi=(X^{1},...,X^{m})\in\mathcal{F}^{m},$ let
\[
\rho_{\xi,x}:\Omega_{\xi,x}\rightarrow\mathbb{R}^{n}:T\mapsto\iota
_{S}\raisebox{2pt}{$\scriptstyle\circ\, $}\xi_{T}(x),
\]
where $\iota_{S}:S\rightarrow\mathbb{R}^{n}$ is the inclusion map. If $M$ is
the orbit of $\mathcal{F}$ through $x$, considered as a subset of
$\mathbb{R}^{n}$, then it is the union of all the images of all the mappings
$\rho_{\xi,x},$ as $m$ varies over the set $\mathbb{N}$ of natural numbers and
$\xi$ varies over $\mathcal{F}^{m}.$ We topologize $M$ by the strongest
topology $\mathcal{T}$ which makes all the maps $\rho_{\xi,x}$ continuous.
Since each $\rho_{\xi,x}:\Omega_{\xi,x}\rightarrow\mathbb{R}^{n}$ is
continuous, it follows that the topology of $M$ as a subspace of
$\mathbb{R}^{n}$ is coarser than the topology $\mathcal{T}$. Hence, the
inclusion of $M$ into $\mathbb{R}^{n}$ is continuous with respect to the
topology $\mathcal{T} $. In particular, $M$ is Hausdorff. Since all the sets
$\Omega_{\xi,x}$ are connected it follows that $M$ is connected. The proof
that the topology $\mathcal{T}$ of $M$ defined above is independent of the
choice of $x\in M$ is exactly the same as in \cite{sussmann}$.$

If $S$ is a subcartesian space, then it can be covered by a family
$\{U_{\alpha}\}_{\alpha\in A}$ of open subsets, each of which is diffeomorphic
to a subset of $\mathbb{R}^{k}$. The argument given above can be repeated in
each $U_{\alpha}$ leading to a topology $\mathcal{T}_{\alpha}$ in $M_{\alpha
}=U_{\alpha}\cap M$. For $\alpha,\beta\in A,$ the topologies $\mathcal{T}%
_{\alpha}$ and $\mathcal{T}_{\beta}$ are the same when restricted to
$M_{\alpha}\cap M_{\beta}$. We define the topology of $M $ so that, for each
$\alpha\in A,$ the induced topology in $M_{\alpha}$ is $\mathcal{T}_{\alpha}.$

Suppose now that $\mathcal{F}$ is a family of vector fields on $S$. For each
$x\in S$, let $D_{\mathcal{F}_{x}}$ be the linear span of $\mathcal{F}%
_{x}=\{X(x)\mid X\in\mathcal{F}\}$. Suppose there is a neighbourhood of $x\in
S$ diffeomorphic to a subset of $\mathbb{R}^{n}$. Then, $\dim D_{\mathcal{F}%
_{x}}\leq n$.

\begin{description}
\item [Lemma 9.]For a locally complete family $\mathcal{F}$ of vector fields
on a subcartesian space $S$, $\dim D_{\mathcal{F}_{x}}$ is constant on orbits
of $\mathcal{F}$.
\end{description}

\noindent\textbf{Proof}. Given $x\in S$, let $\dim D_{\mathcal{F}_{x}}=k$, and
$X^{1},...,X^{k}\in\mathcal{F}$ be such that $\{X^{1}(x),,...,X^{k}(x)\} $ is
a basis in $D_{\mathcal{F}_{x}}.$ Since the family $\mathcal{F}$ is locally
complete, for every $X\in\mathcal{F}$, and $t\in\Omega_{\{X\},x},$
$\varphi_{t\ast}^{X}X^{1}(\varphi_{t}^{X}(x)),...,\varphi_{t\ast}^{X}%
X^{k}(\varphi_{t}^{X}(x))$ are in $D_{\mathcal{F}_{\varphi_{t}^{X}(x)}}$ and
are linearly independent because $\varphi_{t}^{X}$ is a local diffeomorphism.
Hence, $\dim D_{\mathcal{F}_{x}}\leq\dim D_{\mathcal{F}_{\varphi_{t}^{X}(x)}}%
$. Using $\varphi_{-t}^{X}$, we can show that $\dim D_{\mathcal{F}%
_{\varphi_{t}^{X}(x)}}\leq\dim D_{\mathcal{F}_{x}}$. Hence, $\dim
D_{\mathcal{F}_{x}}=\dim D_{\mathcal{F}_{\varphi_{t}^{X}(x)}}$ for every $X\in
F.$ Repeating this argument along $\xi_{T}(x)=\left(  \varphi_{t_{m}}^{X^{m}%
}\raisebox{2pt}{$\scriptstyle\circ
\, $}...\raisebox{2pt}{$\scriptstyle\circ\, $}\varphi_{t_{1}}^{X^{1}}\right)
(x)$, we conclude that $\dim D_{\mathcal{F}_{y}}=\dim D_{\mathcal{F}_{x}}$ for
every $y$ on the orbit of $F$ though $x$. \hfill$\Box\bigskip$

\noindent In analogy with standard terminology, we shall use the term integral
manifold of $D_{\mathcal{F}}$ for a connected manifold $M$ contained in $S,$
such that its inclusion into $S$ is smooth and, for every $x\in M$,
$T_{x}M=D_{\mathcal{F}_{x}}.$

\begin{description}
\item [Theorem 3.]Let $\mathcal{F}$ be a locally complete family of vector
fields on a subcartesian space $S.$ Each orbit $M$ of $\mathcal{F},$ with the
topology $\mathcal{T}$ introduced above, admits a unique manifold structure
such that the inclusion map $\iota_{MS}:M\hookrightarrow S$ is smooth. In
terms of this manifold structure, $M$ is an integral manifold of
$D_{\mathcal{F}}.$
\end{description}

\noindent\textbf{Proof}. Let $M$ be an orbit of $\mathcal{F}$. Since
$\mathcal{F}$ is locally complete, for each $z\in M$, the dimension $m=\dim
D_{\mathcal{F}_{z}}$ is independent of $z$, and there exist $m$ vector fields
$X^{1},...,X^{m}$ in $\mathcal{F}$ that are linearly independent in an open
neighbourhood $V$ of $z$ in $S$. Without loss of generality, we may assume
that $V$ is a subset of $\mathbb{R}^{n}$. By Proposition 2, the restrictions
of $X^{1},...,X^{m}$ to vector fields on $V$ extend to vector fields
$\tilde{X}^{1},...,\tilde{X}^{m}$ on a neighbourhood $U$ of $z$ in
$\mathbb{R}^{n}.$ Without loss of generality, we may assume that they are
linearly independent on $U.$

Given $x\in V\cap U\subseteq\mathbb{R}^{n}$, let $\xi=(X^{1},...,X^{m})$ be
such that $X^{1}(x),...,X^{m}(x)$ form a basis of $D_{\mathcal{F}_{x}},$
$T=(t_{1},...,t_{m})\in\Omega_{\xi,x}$ and
\[
\tilde{\rho}_{\xi,x}:\Omega_{\xi,x}\rightarrow U\subseteq\mathbb{R}%
^{n}:T\mapsto\left(  \varphi_{t_{m}}^{\tilde{X}^{m}}\raisebox{2pt}%
{$\scriptstyle\circ
\, $}...\raisebox{2pt}{$\scriptstyle\circ\, $}\varphi_{t_{1}}^{\tilde{X}^{1}%
}\right)  (x).
\]
For each $i=1,...,m,$
\[
u\frac{d}{dt}\varphi_{t_{i}}^{\tilde{X}^{i}}(x)=u\tilde{X}^{i}(\varphi
_{t}^{\tilde{X}^{i}}(x)).
\]
Hence,
\begin{align*}
&  T\tilde{\rho}_{\xi,x}(T)(u_{1},...,u_{m})=\\
&  =u_{1}\frac{d}{dt_{1}}\left(  \varphi_{t_{m}}^{\tilde{X}^{m}}\raisebox
{2pt}{$\scriptstyle\circ\, $}...\raisebox{2pt}{$\scriptstyle\circ\, $}%
\varphi_{t_{1}}^{\tilde{X}^{1}}\right)  (x)+...+u_{m}\frac{d}{dt_{m}}\left(
\varphi_{t_{m}}^{\tilde{X}^{m}}\raisebox{2pt}{$\scriptstyle\circ
\, $}...\raisebox{2pt}{$\scriptstyle\circ\, $}\varphi_{t_{1}}^{\tilde{X}^{1}%
}\right)  (x)\\
&  =u_{1}(\varphi_{t_{m}}^{\tilde{X}^{m}})_{\ast}...(\varphi_{t_{2}}%
^{\tilde{X}^{2}})_{\ast}\tilde{X}(\varphi_{t_{1}}^{\tilde{X}^{1}%
}(x))+...+u_{m}\tilde{X}^{m}\left(  \varphi_{t_{m}}^{\tilde{X}^{m}}\left(
...\left(  \varphi_{t_{1}}^{\tilde{X}^{1}}(x)\right)  \right)  \right)  .
\end{align*}
In particular,
\begin{equation}
T\tilde{\rho}_{\xi,x}(0)(u_{1},...,u_{m})=u_{1}\tilde{X}^{1}(x)+...+u_{m}%
\tilde{X}^{m}(x). \label{Trho}%
\end{equation}
Since the vectors $\tilde{X}^{1}(x),...,\tilde{X}^{m}(x)$ are linearly
independent, it follows that $T\tilde{\rho}_{\xi,x}(0):\mathbb{R}%
^{m}\rightarrow\mathbb{R}^{n}$ is one to one. Hence, there exists an open
neighbourhood $W_{\xi,x}$ of $0$ in $\mathbb{R}^{m}$ such that the restriction
$\tilde{\rho}_{\xi,x}\mid W_{\xi,x}$ of $\tilde{\rho}_{\xi,x}$ to $W_{\xi,x}$
is an immersion of $W_{\xi,x}$ into $U\subseteq\mathbb{R}^{n}$. Therefore,
$M_{\xi,x}=\tilde{\rho}_{\xi,x}(W_{\xi,x})$ is an immersed submanifold of
$U\subseteq\mathbb{R}^{n}.$ Moreover, there exists a smooth map $\mu_{\xi
,x}:M_{\xi,x}\rightarrow W_{\xi,x}$ such that $\mu_{\xi,x}\raisebox
{2pt}{$\scriptstyle\circ\, $}(\tilde{\rho}_{\xi,x}\mid W_{\xi,x})=identity$.
Every point $y\in M_{\xi,x}$ is of the form
\[
y=\left(  \varphi_{t_{m}}^{\tilde{X}^{m}}\raisebox{2pt}{$\scriptstyle\circ
\, $}...\raisebox{2pt}{$\scriptstyle\circ\, $}\varphi_{t_{1}}^{\tilde{X}^{1}%
}\right)  (x)
\]
for some $T=(t_{1},...,t_{m})\in W_{\xi,x}.$ Since $\tilde{X}^{1}%
,...,\tilde{X}^{m}$ are extensions to $U\subseteq\mathbb{R}^{n}$ of the
restrictions to $V$ of vector fields $X^{1},...,X^{m}$ on $S$, it follows that
$\varphi_{t_{1}}^{\tilde{X}^{1}},...,\varphi_{t_{n}}^{\tilde{X}^{m}}$ preserve
$V$. Moreover, $x\in V$ so that
\[
y=\left(  \varphi_{t_{m}}^{X^{m}}\raisebox{2pt}{$\scriptstyle\circ
\, $}...\raisebox{2pt}{$\scriptstyle\circ\, $}\varphi_{t_{1}}^{X^{1}}\right)
(x)\in V\subseteq S\text{.}%
\]
Hence, $M_{\xi,x}$ is contained in $V\subseteq S$.

Let $\iota_{MV}:M_{\xi,x}\hookrightarrow V$ be the inclusion map. We want to
show that it is smooth. Let $\iota_{M}:M_{\xi,x}\rightarrow U\subseteq
\mathbb{R}^{n}$ and $\iota_{V}:V\rightarrow U\subseteq\mathbb{R}^{n}$ be the
inclusion maps. Then $\iota_{M}=\iota_{V}\raisebox{2pt}{$\scriptstyle\circ
\, $}\iota_{MV}$ and $\iota_{M}^{\ast}f=\iota_{MV}^{\ast}\raisebox
{2pt}{$\scriptstyle\circ\, $}\iota_{V}^{\ast}f$ for every function $f\in
C^{\infty}(U)$. Since $M_{\xi,x}$ is an immersed submanifold of $U,$ it
follows that $\iota_{M}^{\ast}f\in C^{\infty}(M_{\xi,x})$ for all $f\in
C^{\infty}(U)$. Similarly, $V$ is a differential subspace of $U$ so that
$\iota_{V}^{\ast}f\in C^{\infty}(V)$ for all $f\in C^{\infty}(U)$. Moreover,
every $f_{V}\in C^{\infty}(V)$ is locally of the form $\iota_{V}^{\ast}f$ for
some $f\in C^{\infty}(U)$. Since differentiability is a local property, it
follows that $\iota_{MV}^{\ast}f_{V}\in C^{\infty}(M_{\xi,x})$ for every
$f_{V}\in C^{\infty}(V).$ Hence, $\iota_{MV}:M_{\xi,x}\hookrightarrow V$ is smooth.

Thus, for every open set $V$ in $S$, that is diffeomorphic to a subset of
$\mathbb{R}^{n}$, each $x\in V$, and every $\xi=(X^{1},...,X^{m})$ such that
$X^{1}(x),...,X^{m}(x)$ form a basis of $D_{\mathcal{F}_{x}},$ we have a
manifold $M_{\xi,x}$ contained in $V$ such that the inclusion map $\iota
_{MV}:M_{\xi,x}\hookrightarrow V$ is smooth. Since $V$ is open in $S$, the
inclusion of $M_{\xi,x}$ into $S$ is smooth.

Suppose that $M_{\xi_{1},x_{1}}\cap M_{\xi_{2},x_{2}}\ne\emptyset$. If $y\in
M_{\xi_{1},x_{1}}\cap M_{\xi_{2},x_{2}}$, then $y=\tilde{\rho}_{\xi_{1},x_{1}%
}(T_{1})=\tilde{\rho}_{\xi_{2},x_{2}}(T_{2})$ for $T_{1}\in V_{\xi_{1},x_{1}}$
and $T_{2}\in V_{\xi_{2},x_{2}}.$ Hence,
\[
T_{2}=\mu_{\xi_{2},x_{2}}((\tilde{\rho}_{\xi_{1},x_{1}}(T_{1}))).
\]
Since $\tilde{\rho}_{\xi_{1},x_{1}}$ and $\mu_{\xi_{2},x_{2}}$ are smooth, it
follows that the identity map on $M_{\xi_{1},x_{1}}\cap M_{\xi_{2},x_{2}}$ is
a diffeomorphism of the differential structures on $M_{\xi_{1},x_{1}}\cap
M_{\xi_{2},x_{2}}$ induced by the inclusions into $M_{\xi_{1},x_{1}}$ and
$M_{\xi_{2},x_{2}}$, respectively. Therefore $M_{\xi_{1},x_{1}}\cup M_{\xi
_{2},x_{2}}$ is a manifold contained in $S$ and the inclusion of $M_{\xi
_{1},x_{1}}\cup M_{\xi_{2},x_{2}}$ into $S$ is smooth.

Since $M=\bigcup_{\xi,x}M_{\xi,x},$ the above argument shows that $M$ is a
manifold contained in $S$ such that the inclusion map $M\hookrightarrow S$ is
smooth. Moreover, the manifold topology of $M$ agrees with the topology
$\mathcal{T}$ discussed above. Finally, equation (\ref{Trho}) implies that $M$
is an integral manifold of $D_{\mathcal{F}}.$ \hfill$\Box$

\bigskip

We see from Theorem 3 that a locally complete family $\mathcal{F}$ of vector
fields on a subcartesian space $S$ gives rise to a partition of $S$ by orbits
of $\mathcal{F}.$ We shall refer to such a partition as a singular foliation
of $S$.

\begin{description}
\item [Theorem 4.]Orbits of the family $\mathcal{X}(S)$ of all vector fields
on a subcartesian space $S$ are manifolds. For every family $\mathcal{F}$ of
vector fields on $S$, orbits of $\mathcal{F}$ are contained in orbits of
$\mathcal{X}(S).$
\end{description}

\noindent\textbf{Proof}. We have shown in Theorem 2 that the family
$\mathcal{X}(S)$ of all vector fields on $S$ is locally complete. Hence, it
gives rise to a partition of $S$ by manifolds. If $\mathcal{F}$ is a family of
vector fields on $S,$ then $\mathcal{F}\subseteq\mathcal{X}(S)$, and every
orbit of $\mathcal{F}$ is contained in an orbit of $\mathcal{X}(S).$
\hfill$\Box$

\bigskip

Theorem 4 asserts that the singular foliation of a subcartesian space $S$ by
orbits of the family $\mathcal{X}(S)$ of all vector fields on $S$ is coarsest
within the class of singular foliations given by orbits of locally complete
families of vector fields. The following example shows that there may be
partitions of a differential space into manifolds which are coarser than the
singular foliation by orbits of $\mathcal{X}(S).$

\begin{description}
\item [Example 3.]Let $S=M_{1}\cup M_{2}\cup M_{3}\cup M_{4}$, where
\begin{align*}
M_{1}  &  =\{(x,y)\in\mathbb{R}^{2}\mid x>0\text{ and }y=\sin(x^{-1})\},\\
M_{2}  &  =\{(x,y)\in\mathbb{R}^{2}\mid x=0\text{ and }-1<y<1\},
\end{align*}
$M_{3}=\{(0,-1)\}$ and $M_{4}=\{(0,1)\}$. Clearly, $M_{1}$ and $M_{2}$ are
manifolds of dimension 1, while $M_{3}$ and $M_{4}$ are manifolds of dimension
0. However, for every $(0,y)\in M_{2}$, a vector $u\in T_{(0,y)}M_{2}$ can be
extended to a vector field on $S$ only if $u=0$. Hence, orbits of the minimal
singular foliation of $S$ are $M_{1}$ and singletons $\{(0,y)\}$ for $-1\leq
y\leq1.$
\end{description}

Having established the existence of the singular foliation of $S$ by orbits of
$\mathcal{X}(S)$, we can study arbitrary families of vector fields.

\begin{description}
\item [Theorem 5.]Let $\mathcal{F}$ be a family of vector fields on a
subcartesian space $S.$ For every $x\in S,$ the orbit $N$ of $F$ through $x$
is a manifold such that the inclusion map $\iota_{NS}:N\hookrightarrow S$ is smooth.
\end{description}

\noindent\textbf{Proof}. Since $\mathcal{F}\subseteq\mathcal{X}(S),$ the
accessible set $N$ of $D_{\mathcal{F}}$ through $x$ is contained in the orbit
$M$ of $\mathcal{X}(S)$ through $x.$ Let $D_{\mathcal{F}}\mid M$ be the
restriction of $D_{\mathcal{F}}$ to the points of $M$. It is a generalized
distribution on $M$, and $N$ is an accessible set of $D_{\mathcal{F}%
}\mathcal{\mid}M.$ It follows from Sussmann's theorem, \cite{sussmann}, that
$N$ is a manifold and the inclusion map $\iota_{NM}:N\rightarrow M$ is smooth.
Since the inclusion $\iota_{MS}:M\rightarrow S$ is smooth, it follows that
$\iota_{NS}=\iota_{MS}\raisebox{2pt}{$\scriptstyle\circ\, $}\iota
_{NM}:N\hookrightarrow S$ is smooth. \hfill$\Box$

\section{Stratified spaces}

In this section we show that smooth stratifications are subcartesian spaces.
This enables us to use the results of the preceding sections in discussing
stratified spaces. For a comprehensive study of stratified spaces see
\cite{goresky-macpherson}, \cite{pflaum} and the references quoted there.

Let $S$ be a paracompact Hausdorff space. A stratification of $S$ is given by
a locally finite partition of $S$ into locally closed subspaces $M\subseteq
S,$ called strata, satisfying the following conditions:

\begin{description}
\item [Manifold Condition.]Every stratum $M$ of $S$ is a smooth manifold in
the induced topology.

\item[Frontier Condition.] If $M$ and $N$ are strata of $S$ such that the
closure $\bar{N}$ of $N$ has a non-empty intersection with $M$, then
$M\subset\bar{N}.$
\end{description}

A smooth chart\textit{\ }on a stratified space $S$ is a homeomorphism
$\varphi$ of an open set $U\subseteq S$ to a subspace $\varphi(U)$ of
$\mathbb{R}^{n}$ such that, for every stratum $M$ of $S$, the image
$\varphi(U\cap M)$ is a smooth submanifold of $\mathbb{R}^{n}$ and the
restriction $\varphi\mid U\cap M:U\cap M\rightarrow\varphi(U\cap M)$ is
smooth. As in the case of manifolds, one introduces the notion of
compatibility of smooth charts, and the notion of a maximal atlas of
compatible smooth charts on $S$. A smooth structure on $S$ is given by a
maximal atlas of smooth charts on $S.$ A continuous function $f:S\rightarrow
\mathbb{R}$ is smooth if, for every $x\in S$ and every chart $\varphi
:U\rightarrow\mathbb{R}^{n}$ with $x\in U,$ there exists an neighbourhood
$U_{x}$ of $x$ contained in $U$ and a smooth function $g:\mathbb{R}%
^{n}\rightarrow\mathbb{R}$ such that $f\mid U_{x}=g\raisebox{2pt}%
{$\scriptstyle\circ\, $}\varphi\mid U_{x}.$ For details see ( \cite{pflaum},
sec. 1.3).

Stratifications can be ordered by inclusion. If we have two stratifications of
the same space $S$, we say that the first stratification is smaller than the
second if every stratum of the second stratification is contained in a stratum
of the first one. For a stratified space $S$, there exists a minimal
stratification of $S$. Some authors reserve the term stratification for a
minimal stratification.

\begin{description}
\item [Theorem 6.]A smooth stratified space is a subcartesian space.
\end{description}

\noindent\textbf{Proof}. Let $S$ be a smooth stratified space and $C^{\infty
}(S)$ the space of smooth functions on $S$ defined above. First, we need to
show that the family $C^{\infty}(S)$ satisfies the conditions given at the
beginning of section 2.

A family $\{W_{\alpha}\}_{\alpha\in A}$ of open sets on $S$ is a subbasis for
the topology of $S$ if, for each $x\in S$ and each open neighbourhood $V$ of
$x$ in $S$, there exist $\alpha_{1},...,\alpha_{p}\in A$ such that $x\in
W_{\alpha_{1}}\cap...\cap W_{\alpha_{p}}\subseteq V.$ Given $x\in S,$ there
exists a chart $\varphi$ on $S$ with domain $U$ containing $x$. If $V$ is a
neighbourhood of $x$ in $S$, then the restriction of $\varphi$ to $V\cap U$ is
a homeomorphism on a set $\varphi(V\cap U)$ in $\mathbb{R}^{n}$ containing
$\varphi(x).$ There exists an open neighbourhood $W$ of $x$ in $V\cap U$ such
that $\varphi(x)\in\varphi(W)\subseteq\overline{\varphi(W)}\subseteq
\varphi(V\cap U)$ and a function $f\in C^{\infty}(S)$ such that $f\mid W=1$
and $f\mid S\backslash(V\cap U)=0$. Hence, $x\in f^{-1}((0,2))\subseteq V$.
This implies that condition 2.1 is satisfied.

Suppose that $f_{1},...,f_{n}\in C^{\infty}(S)$ and $F:C^{\infty}%
(\mathbb{R}^{m})$. We want to show that $F(f_{1},...,f_{m})\in C^{\infty}(S)$.
For every $x\in S$ and every chart $\varphi:U\rightarrow\mathbb{R}^{n}$ with
$x\in U,$ there exists a neighbourhood $U_{x}$ of $x$ contained in $U$ and
smooth functions $g_{1},...,g_{m}:\mathbb{R}^{n}\rightarrow\mathbb{R}$ such
that $f_{1}\mid U_{x}=g_{i}\raisebox{2pt}{$\scriptstyle\circ\, $}\varphi\mid
U_{x}$ for $i=1,...,m$. Hence,
\begin{align*}
F(f_{1},...f_{m})  &  \mid U_{x}=F(f_{1}\mid U_{x},...,f_{m}\mid
U_{x})=F(g_{1}\raisebox{2pt}{$\scriptstyle\circ\, $}\varphi\mid U_{x}%
,...,g_{m}\raisebox{2pt}{$\scriptstyle\circ\, $}\varphi\mid U_{x})\\
&  =F(g_{1},...,g_{n})\mid\raisebox{2pt}{$\scriptstyle\circ\, $}\varphi\mid
U_{x}%
\end{align*}
and condition 2.2 is satisfied.

In order to prove condition 2.3, consider $f:S\rightarrow\mathbb{R}$ such
that, for every $x\in S$, there exists an open neighbourhood $W_{x}$ of $x$
and a function $f_{x}\in C^{\infty}(S)$ satisfying
\begin{equation}
f_{x}\mid W_{x}=f\mid W_{x}. \label{germ}%
\end{equation}
Given $x\in S$, let $f_{x}\in C^{\infty}(S)$ and an open neighbourhood $W_{x}$
be such that equation (\ref{germ}) is satisfied. Let $\varphi:U\rightarrow
\mathbb{R}^{n}$ be a chart such that $x\in U$. There exists an open
neighbourhood $U_{x}$ of $x$ contained in $W_{x}\cap U$ and a smooth function
$g_{x}:\mathbb{R}^{n}\rightarrow\mathbb{R}$ such that $f_{x}\mid U_{x}%
=g_{x}\raisebox{2pt}{$\scriptstyle\circ\, $}\varphi\mid U_{x}$. Since
$U_{x}\subseteq W_{x}\cap U$, it follows from equation (\ref{germ}) that
$f\mid U_{x}=g_{x}\raisebox{2pt}{$\scriptstyle\circ\, $}\varphi\mid U_{x}.$
This holds for every $x\in S,$ which implies that $f\in C^{\infty}(S).$

We have shown that smooth functions on $S$ satisfy the conditions for a
differential structure on $S$. Thus, $S$ is a differential space. Local charts
are local diffeomorphisms of $S$ onto subsets of $\mathbb{R}^{n}$. This
implies that $S$ is a subcartesian space. \hfill$\Box$

A stratified space $S$ is said to be topologically locally trivial\textit{\ }%
if, for every $x\in S$, there exists an open neighbourhood $U$ of $x$ in $S,$
a stratified space $F$ with a distinguished point $o\in F$ such that the
singleton $\{o\}$ is a stratum of $F,$ and a homeomorphism $\varphi
:U\rightarrow(M\cap U)\times F$, where $M$ is the stratum of $S$ containing
$x$, such that $\varphi$ induces smooth diffeomorphisms of the corresponding
strata, and $\varphi(y)=(y,o)$ for every $y\in M\cap U$. The stratified space
$F$ is called the typical fibre over $x$. Since we are dealing here with the
$C^{\infty}$ category, we shall say that a smooth stratified space $S$ is
locally trivial if it is topologically locally trivial and, for each $x\in S$,
the typical fibre $F$ over $S$ is smooth and the homeomorphism $\varphi
:U\rightarrow(M\cap U)\times F$ is a diffeomorphism of differential spaces. In
\cite{cushman-sniatycki} we have shown that the orbit space of a proper action
is locally trivial.

The stratified tangent bundle $T^{s}S$ of a stratified space $S$ is the union
of tangent bundle spaces $TM$ of all strata $M$ of $S.$ We denote by
$\tau:T^{s}S\rightarrow S$ the projection map such that for every $x\in S,$
$\tau^{-1}(x)=T_{x}M$, where $M$ is the stratum containing $x$. For each chart
$\varphi$ on $S,$ with domain $U$ and range $V\subseteq\mathbb{R}^{n},$ one
sets $T^{s}U=\tau^{-1}(U)$ and defines $T\varphi:T^{s}U\rightarrow
T^{s}V\subseteq\mathbb{R}^{2n}$ by requiring that $(T\varphi)\mid TM\cap
T^{s}U=T(\varphi\mid M\cap U)$ for all strata $M$ of $S$. One supplies
$T^{s}S$ with the coarsest topology such that all $T^{s}U\subseteq T^{s}S$ are
open and all $T\varphi$ are continuous, see \cite{pflaum}. A stratified vector
field on $S$ is a continuous section $X$ of $\tau$ such that, for every
stratum $M$ of $S$, the restriction $X\mid M$ is a smooth vector field on $M$.

Let $S$ be a smooth stratified space. By Theorem 6, it is a subcartesian
space. The above definition of a stratified vector field does not ensure that
it generates local one-parameter groups of local diffeomorphisms of $S$.
Conversely, one often uses the term stratification for a partition of a smooth
manifold $S$ which satisfies the Manifold Condition and the Frontier
Condition. In this case, there exist vector fields on $S$, that generate local
one-parameter groups of local diffeomorphisms of $S$, but are not stratified
in the sense given above. In this paper, we shall use the term strongly
stratified vector field on $S$ for a vector field $X$ on $S$ that generates a
local one-parameter group of local diffeomorphisms of $S$, and is such that,
for every stratum $M$ of $S$, $X$ restricts to a smooth vector field on $M$.
Thus, a strongly stratified vector field on $S$ generates a local
one-parameter group of local diffeomorphisms of $S$ that preserves the
stratification structure of $S$.

\begin{description}
\item [Lemma 10.]Let $S$ be a locally trivial stratified space, and $X_{M}$ a
smooth vector field on a stratum $M$ of $S$. For every $x\in M$, there exists
a neighbourhood $W$ of $x$ in $M$ and a strongly stratified vector field $X$
on $S$ such that $X_{M}\mid W=X\mid W.$
\end{description}

\noindent\textbf{Proof.} Let $\varphi_{t}$ be the local one-parameter group of
local diffeomorphisms of $M$ generated by $X_{M}$. Since $S$ is locally
trivial, there exists a neighbourhood $U$ of $x$ in $S,$ a smooth stratified
space $F$ and a diffeomorphism $\varphi:U\rightarrow(M\cap U)\times F$ such
that $\varphi$ induces smooth diffeomorphisms of the corresponding strata, and
$\varphi(y)=(y,o)$ for every $y\in M\cap U$. Each stratum of $(M\cap U)\times
F$ is of the form $(M\cap U)\times N$, where $N$ is a stratum of $F$. Let
$X_{U}$ be a stratified vector field on $(M\cap U)\times F$ such that, for
every stratum $N$ of $F$,
\[
X_{U}\mid(M\cap U)\times N=(X_{M}\mid M\cap U,0).
\]
It is a derivation generating a local one-parameter group of local
diffeomorphisms $\psi_{t}$ of $(M\cap U)\times F$ such that, for every
$(y,z)\in(M\cap U)\times F,$ $\psi_{t}(y,z)$ is defined whenever $\varphi
_{t}(y)$ is defined, is contained in $M\cap U$, and $\psi_{t}(y,z)=(\varphi
_{t}(y),z)$. Hence, $X_{U}$ is a vector field on $(M\cap U)\times F$.

Let $V_{1}$ and $V_{2}$ be neighbourhoods of $x$ in $S$ such that $\bar{V}%
_{1}\subseteq V_{2}\subseteq\bar{V}_{2}\subseteq U$. There exists a function
$f\in C^{\infty}(S)$ such that $f\mid\bar{V}_{1}=1$ and $f\mid S\backslash
V_{2}=0$. Let $X$ be a stratified vector field on $S$ such that
\begin{align*}
X(x^{\prime})  &  =f(x^{\prime})(T\varphi^{-1}(X_{U}(\varphi(x^{\prime
}))))\text{ for }x^{\prime}\in U,\\
X(x^{\prime})  &  =0\text{ for }x^{\prime}\in S\backslash U.
\end{align*}
Since $X_{U}$ is a vector field on $(M\cap U)\times F,$ it follows that $X$ is
a vector field on $S$. Let $W=V_{1}\cap M$. Since $f\mid W=1,$ it follows that
$X\mid W=X_{M}\mid W$, which completes the proof. \hfill$\Box$

\bigskip

Let $\mathcal{X}_{s}(S)$ denote the family of all strongly stratified vector
fields on a smooth stratified space $S$.

\begin{description}
\item [Lemma 11.]The family $\mathcal{X}_{s}(S)$ of all strongly stratified
vector fields on a locally trivial stratified space $S$ is locally complete.
\end{description}

\noindent\textbf{Proof}. For $X\in\mathcal{X}_{s}(S)$, let $\varphi_{t}$
denote the local one-parameter group of local diffeomorphisms of $S$ generated
by $X$. Suppose $U$ is the domain of $\varphi_{t}$ and $V$ is its range. In
other words $\varphi_{t}$ maps $U\,$ diffeomorphically onto $V$. In Lemma 8 we
have shown that, for each $Y\in\mathcal{X}_{s}(S)$, $\varphi_{t\ast}Y$ is a
vector field on $V$. By Lemma 1, for every $x\in V$, there exists an open
neighbourhood $W$ of $x$ in $V$ such that $\bar{W}\subseteq V$ and a function
$f\in C^{\infty}(S)$ such that $f\mid W=1$ and $f\mid S\backslash\bar{V}=0$.
Hence, there is a vector field $Z$ on $S$ such that $Z\mid V=(f\mid
V)\varphi_{t\ast}Y$ and $Z\mid S\backslash\bar{V}=0$. In particular, $Z\mid
W=\varphi_{t\ast}Y\mid W.$

For every stratum $M$ of $S$, the restriction of $Y$ to $M$ is tangent to $M$.
Moreover, $X\in\mathcal{X}_{c}(S)$ implies that $\varphi_{t*}$ preserves $M$.
Hence, $\varphi_{t*}Y$ restricted to $V\cap M$ is tangent to $V\cap M$. Since
$Z\mid V=\varphi_{t*}Y\mid V,$ it follows that $Z\mid V\cap M$ is tangent to
$V\cap M$. On the other hand $Z\mid S\backslash\bar{V}=0,$ which implies that
$Z\mid(S\backslash\bar{V})\cap M=0$ is tangent to $(S\backslash\bar{V})\cap
M$. Hence, $Z\mid M$ is tangent to $M.$ This ensures that $Z$ is a strongly
stratified vector field on $S$.

The argument above is valid for every $X$ and $Y$ in $\mathcal{X}_{s}(S)$.
Hence, $\mathcal{X}_{s}(S)$ is a locally complete family of vector fields on
$S$. \hfill$\Box$

\begin{description}
\item [Theorem 7.]Strata of a locally trivial stratified space $S$ are orbits
of the family $\mathcal{X}_{s}(S)$ of strongly stratified vector fields.
\end{description}

\noindent\textbf{Proof}. By Lemma 11, the family $\mathcal{X}_{s}(S)$ of all
strongly stratified vector fields on $S$ is complete. Hence, its orbits give
rise to a singular foliation of $S$. By definition, for each stratum $M$ of
$S$ and every $X\in\mathcal{X}_{s}(S),$ the restriction of $X$ to $M$ is
tangent to $M$. Hence, orbits of $X_{s}(S)$ are contained in strata of $S$.

Let $x$ and $y$ be in the same stratum $M$ of $S$. Since $M$ is connected,
there exists a piecewise smooth curve $c$ in $M$ joining $x$ to $y$. In other
words, there exist vector fields $X_{M}^{1},...,X_{M}^{l}$ on $M$ such that
$y=\varphi_{t_{1}}^{l}\raisebox{2pt}{$\scriptstyle\circ\, $}...\raisebox
{2pt}{$\scriptstyle\circ
\, $}\varphi_{t_{l}}^{1}(x)$, where $\varphi_{t}^{i}$ is the local
one-parameter group of local diffeomorphisms of $M$ generated by $X_{M}^{i}$
for $i=1,...,l.$ Let $x_{t}^{i}=\varphi_{t}^{i}\raisebox{2pt}{$\scriptstyle
\circ
\, $}\varphi_{t_{i-1}}^{i-1}\raisebox{2pt}{$\scriptstyle\circ\, $}%
...\raisebox{2pt}{$\scriptstyle\circ\, $}\varphi_{t_{1}}^{1}$ for every
$i=1,...,l$ and $t\in\lbrack0,t_{i}]$. There exists a neighbourhood $W_{t}%
^{i}$ of $x_{t}^{i}$ in $M$ and a vector field $X_{t}^{i}$ on $S$ such that
$X_{t}^{i}\mid W_{t}^{i}=X_{M}^{i}\mid W_{t}^{i}$. The family $\{W_{t}^{i}\mid
i=1,...l,$ $t\in\lbrack0,t_{i}]\}$ gives a covering of the curve $c$ joining
$x$ to $y$ $.$Since the range of $c$ is compact, there exists a finite
subcovering $\{W_{t_{j}}^{i}\mid i=1,...l,$ $j=1,...,n_{i}\}$ covering $c$.
Hence, $c$ is a piecewise integral curve of the vector fields $X_{t_{j}}^{i},$
$i=1,...l,$ $j=1,...,n_{i}$, on $S$. This implies that $M$ is an
$\mathcal{X}(S)$ orbit. \hfill$\Box$

\begin{description}
\item [Theorem 8.]Let $S$ be a locally trivial stratified space. Then the
singular foliation of $S$ by orbits of the family $\mathcal{X}(S)$ of all
vector fields on $S$ is a smooth stratification.
\end{description}

\noindent\textbf{Proof}. By Theorem 7, strata of the stratification of $S$ are
orbits of a family of vector fields on $S$. Hence, each stratum of the
stratification of $S$ is contained in an orbit of the family $\mathcal{X}(S)$
of all vector fields on $S$. Thus, every orbit of $X(S)$ is the union of
strata. Since strata of $S$ form a locally finite partition of $S$, it follows
that the singular foliation of $S$ by orbits of $\mathcal{X}(S)$ is also
locally finite.

Next, we show that orbits of $\mathcal{X}(S)$ are locally closed. Let $P$ be
an orbit of $X(S)$ through $x\in S,$ and $M$ the stratum of $S$ containing $x
$. Let $W_{0}$ be a neighbourhood of $x$ in $S$ which intersects a finite
number of strata $M_{1},...,M_{n}$ of $S$. In other words,
\[
W_{0}=\bigcup_{i=1}^{n}W_{0}\cap M_{i}.
\]
Then
\[
U=W_{0}\cap P=\bigcup_{i=1}^{n}W_{0}\cap M_{i}\cap P=\bigcup_{M_{i}\subseteq
P}W_{0}\cap M_{i}%
\]
is a neigbourhood of $x$ in $P$. Each $W_{0}\cap M_{i}$ is an open subset of
$M_{i}$. Since $M_{i}$ is locally closed, we can choose $W_{0}$ sufficiently
small so that there exists a closed set $V_{i}$ in $S$ such that $W_{0}\cap
M_{i}=W_{0}\cap V_{i}$. Hence,
\[
U=\bigcup_{M_{i}\subseteq P}W_{0}\cap M_{i}=\bigcup_{M_{i}\subseteq P}%
W_{0}\cap V_{i}=W_{0}\cap\left(  \bigcup_{M_{i}\subseteq P}V_{i}\right)  ,
\]
and
\[
V=\bigcup_{M_{i}\subseteq P}V_{i}%
\]
is closed as a finite union of closed sets. This shows that $P$ is locally closed.

It remains to verify the Frontier Condition. Let $P$ and $Q$ be orbits of
$\mathcal{X}(S)$ such that there exists a point $x\in\bar{P}\cap Q$, where
$\bar{P}$ is the closure of $P$. We want to show that $Q\subseteq\bar{P}$.
Suppose that $Q$ is not contained in $\bar{P}$. Then there exists a point
$y\in Q$ such that $y\notin\bar{P}$. Since $Q$ is an orbit of $\mathcal{X}(S)
$, there exists a piecewise smooth curve $\gamma:[0,1]\rightarrow Q$ such that
$x=\gamma(0)$ and $y=\gamma(1)$. Let
\[
\tau=\sup\{t\in\lbrack0,1]\mid\gamma(s)\in\bar{P}\text{ for all }s\leq t\}.
\]
Since $y=\gamma(1)\notin\bar{P}$, it follows that $0\leq\tau<1.$ Then
$z=\gamma(\tau)\in\bar{P}$ and, for every $\varepsilon>0$, there exists
$t\in(0,\varepsilon)$ such that $\gamma(\tau+t)\notin\bar{P}.$ Since $\gamma$
is piecewise smooth, there exists a vector field $X$ on $S$ such that
\[
X(z)=\lim_{t\rightarrow0+}\dot{\gamma}(\tau+t)\neq0.
\]
Let $\varphi_{t}^{X}$ denote the local one-parameter group of local
diffeomorphisms of $S$ generated by $X$. For sufficiently small $t_{0}>0,$ we
can choose a neighbourhood $U$ of $z$ in $S$ such that $\varphi_{t}^{X}$ is
defined on $U$ for all $0\leq t\leq t_{0}$. Let $t\in(0,t_{0})$ be such that
$\gamma(\tau+t)\notin\bar{P}$. We have
\[
\gamma(\tau+t)=\varphi_{t}^{X}(\gamma(\tau))=\varphi_{t}^{X}(z).
\]
Since $z\in\bar{P}$, it follows that $U\cap P\neq\emptyset$. Since $P$ is an
orbit of $\mathcal{X}(S),$ it follows that $\varphi_{s}^{X}(U\cap P)\subseteq
P$ for all $0\leq s\leq t.$ Moreover, $\varphi_{t}^{X}$ is a diffeomorphism of
$U$ on its image mapping $U\cap P$ onto $\varphi_{t}^{X}(U\cap P)$. Hence
\[
\gamma(\tau+t)=\varphi_{t}^{X}(z)\in\varphi_{t}^{X}(\bar{P}\cap U)\subseteq
\overline{\varphi_{t}^{X}(U\cap P)}\subseteq\bar{P},
\]
which contradicts the assumption that $\gamma(\tau+t)\notin\bar{P}$. This
implies that $Q\subseteq\bar{P},$ which completes the proof that the singular
foliation of $S$ by orbits of the family $\mathcal{X}(S)$ of all vector fields
on $S$ is a stratification.

We still need to show that the stratification of $S$ by orbits of the family
$\mathcal{X}(S)$ of all vector fields on $S$ is smooth. By assumption, $S$ is
a smooth stratified space, hence its smooth structure is determined by a
maximal atlas of compatible smooth charts on $S$. Let $\varphi$ be a chart of
this atlas with domain $U$ and range $\varphi(U)\subseteq\mathbb{R}^{n}$. For
each orbit $M$ of $\mathcal{X}(S),$ the intersection $M\cap U$ is a manifold
contained in $U$, and $\varphi(M\cap U)$ is a manifold contained in
$\varphi(U)$ because $\varphi$ is a diffeomorphism. Suppose that
$f:\varphi(U\cap M)\rightarrow\mathbb{R}$ is smooth. Since $M$ is locally
closed, it follows that $\varphi(U\cap M)$ is a locally closed subset of
$\mathbb{R}^{n}$. Hence, for every $y\in\varphi(U\cap M)$, there exists a
neighbourhood $V$ of $y$ in $\mathbb{R}^{n}$ and $\tilde{f}\in C^{\infty
}(\mathbb{R}^{n})$ such that $f\mid V\cap\varphi(U\cap M)=\tilde{f}\mid
V\cap(\varphi(U\cap M).$ Since $\varphi$ is a smooth map of $U$ into
$\mathbb{R}^{n}$, $\varphi^{\ast}\tilde{f}=\tilde{f}\raisebox{2pt}%
{$\scriptstyle\circ\, $}\varphi\in C^{\infty}(U)$. Moreover,
\begin{align*}
(\varphi &  \mid U\cap M)^{\ast}f\mid\varphi^{-1}(V)=f\raisebox{2pt}%
{$\scriptstyle\circ\, $}\varphi\mid\varphi^{-1}(V\cap\varphi(U\cap M))=f\mid
V\cap\varphi(U\cap M)\\
&  =\tilde{f}\mid V\cap(\varphi(U\cap M)=\tilde{f}\raisebox{2pt}{$\scriptstyle
\circ\, $}\varphi\mid\varphi^{-1}(V\cap\varphi(U\cap M))=\tilde{f}%
\raisebox{2pt}{$\scriptstyle\circ\, $}\varphi\mid\varphi^{-1}(V)\cap(U\cap M).
\end{align*}
Since this is valid for every $y\in\varphi(U\cap M),$ it follows that
$\varphi\mid U\cap M:U\cap M\rightarrow$ $\varphi(U\cap M)$ is smooth. This
holds for every chart $\varphi$ and every orbit $M$ of $\mathcal{X}(S)$.
Hence, the stratification of $S$ by orbits of $\mathcal{X}(S)$ is smooth with
respect to the atlas defining the smooth structure on $S$. \hfill$\Box$

\begin{description}
\item [Remark]In Theorem 8, we have assumed that the original stratification
of $S$ is locally trivial. We do not know if the stratification of $S$ given
by orbits of $\mathcal{X}(S)$ is locally trivial.
\end{description}

\section{Poisson reduction}

We can now return to the problem of Poisson reduction of a proper symplectic
action
\[
\Phi:G\times P\rightarrow P:(g,p)\mapsto\Phi(g,p)\equiv\Phi_{g}(p)\equiv gp
\]
of a Lie group $G$ on a symplectic manifold $(P,\omega),$ which has motivated
this work. Here, $\omega$ is a symplectic form on $P$ and, for every $g\in G$,
$\Phi_{g}^{\ast}\omega=\omega$.

For every $p\in P$, the isotropy group $G_{p}$ of $p$ is given by
\[
G_{p}=\{g\in G\mid\Phi_{g}(p)=p\}.
\]
Since the action $\Phi$ is proper, all isotropy groups are compact. For every
compact subgroup $K$ of $G$,the set
\[
P_{K}=\{p\in P\mid G_{p}=K\}
\]
of points of isotropy type $K$, and the set
\[
P_{(K)}=\{p\in P\mid G_{p}\text{ is conjugate to }K\}
\]
of points of orbit type $K$ are local manifolds. Thus means that connected
components of $P_{K}$ and $P_{(K)}$ are submanifolds of $P,$
\cite{bates-cushman}.

Let $S=P/G$ denote the orbit space of the action $\Phi$ and $\rho:P\rightarrow
S$ the orbit the orbit map associating, to each $p\in P,$ the orbit
$Gp=\{\Phi_{g}(p)\mid g\in G\}$ of $G$ through $p.$ The orbit space $S $ is
stratified by orbit type, \cite{D-K}. In other words, strata of $S$ are
connected components of $\rho(P_{(K)}),$ as $K$ varies over compact subgroups
of $G$ for which $P_{(K)}\neq\emptyset.$ The space $S$ is a differential space
with a differential structure $C^{\infty}(S),$ introduced in \cite{schwarz},
which consists of push-forwards to $S$ of $G$-invariant smooth functions on
$P$, it is locally trivial and minimal, \cite{bierstone1}, \cite{bierstone2}.
Minimality of the stratification implies that its strata are orbits of the
family $\mathcal{X}(S)$ of all vector fields on $S.$

For each $f\in C^{\infty}(P)$, the Hamiltonian vector field of $f$ is the
vector field $X_{f}$ on $P$ defined by
\begin{equation}
X_{f}%
{\mbox{$ \rule{5pt} {.5pt}\rule{.5pt} {6pt}\, $}}%
\omega=df, \label{Hamiltonian0}%
\end{equation}
where $%
{\mbox{$ \rule{5pt} {.5pt}\rule{.5pt} {6pt}\, $}}%
$ denotes the left interior product of vector fields and forms. The Poisson
bracket of $f_{1},f_{2}\in C^{\infty}(P)$ is given by
\[
\{f_{1},f_{2}\}=X_{f_{1}}\cdot f_{2}.
\]
It is antisymmetric, satisfies the Jacobi identity
\[
\{f_{1},\{f_{2},f_{3}\}\}+\{f_{2},\{f_{3},f_{1}\}\}+\{f_{3},\{f_{1}%
,f_{2}\}\}=0,
\]
and the Leibniz rule
\[
\{\{f_{1},f_{2}f_{3}\}=\{f_{1},f_{2}\}f_{3}+f_{2}\{f_{1},f_{3}\},
\]
for every $f_{1},f_{2},f_{3}\in C^{\infty}(P).$

Since the action $\Phi$ of $G$ on $P$ preserves the symplectic form $\omega$,
the induced action of $G$ on $C^{\infty}(P)$ preserves the Poisson bracket.
Hence, the space $C^{\infty}(P)^{G}$ of $G$-invariant functions in $C^{\infty
}(P)$ is a Poisson sub-algebra of $C^{\infty}(P)$. The differential structure
\[
C^{\infty}(S)=\{h:S\rightarrow\mathbb{R}\mid\rho^{\ast}h\in C^{\infty}(P)\}
\]
is isomorphic to $C^{\infty}(P)^{G}.$ This implies that the Poisson bracket on
$C^{\infty}(P)^{G}$ induces a Poisson bracket on $C^{\infty}(S)$ such that
\begin{equation}
\rho^{\ast}\{h_{1},h_{2}\}=\{\rho^{\ast}h_{1},\rho^{\ast}h_{2}\}
\label{Poisson}%
\end{equation}
for every $h_{1},h_{2}\in C^{\infty}(S).$

For every $h\in C^{\infty}(S)$, we denote by $X_{h}$ the derivation of
$C^{\infty}(S)$ given by
\[
X_{h}\cdot f=\{h,f\}\text{ for all }f\in C^{\infty}(S).
\]
Sjamaar and Lerman showed that, for every $x\in S$, there exists a unique
maximal integral curve $\gamma$ of $X_{h}$ through $x$, \cite{sjamaar-lerman}.
For every $f\in C^{\infty}(S)$,
\[
\rho^{\ast}(X_{h}\cdot f)=X_{\rho^{\ast}h}\cdot\rho^{\ast}f,
\]
where $X_{\rho^{\ast}h}$ is the Hamiltonian vector field of $\rho^{\ast}h\in
C^{\infty}(P)$ defined by equation (\ref{Hamiltonian0}). Since $\rho^{\ast}h$
is $G$-invariant, the one-parameter local group $\varphi_{t}^{X_{\rho^{\ast}%
h}}$ of local diffeomorphisms of $P$ generated by $X_{\rho^{\ast}h}$ commutes
with the action of $G.$ If $p\in\rho^{-1}(x)\subseteq P$, then $\gamma
(t)=\rho\raisebox{2pt}{$\scriptstyle\circ\, $}\varphi_{t}^{X_{\rho^{\ast}h}%
}(x).$ Hence, translations along integral curves of $X_{h}$ give rise to a
local one-parameter group $\varphi_{t}^{X_{h}}$ of local diffeomorphisms of
$S$ such that $\rho\raisebox{2pt}{$\scriptstyle\circ\, $}\varphi_{t}%
^{X_{\rho^{\ast}h}}=\varphi_{t}^{X_{h}}\raisebox{2pt}{$\scriptstyle\circ
\, $}\rho$. This implies that the derivation $X_{h}$ of $C^{\infty}(S)$ is a
vector field on $S$ in the sense of the definition adopted in Section 4. We
shall refer to $X_{h}$ as the Hamiltonian vector field on $S$ corresponding to
$h\in C^{\infty}(S).$

Hamiltonian vector fields on $(P,\omega)$ preserve the symplectic form
$\omega$. Hence, they preserve the Poisson bracket on $(P,\omega)$. In other
words,
\[
(\varphi_{-t}^{X_{h}})^{\ast}\{(\varphi_{t}^{X_{h}})f_{1},(\varphi_{t}^{X_{h}%
})f_{2}\}=\{f_{1},f_{2}\}\text{ for all }f_{1},f_{2},h\in C^{\infty}(P).
\]
Restricting this equality to $G$-invariant functions, and taking into account
the definition of the Poisson bracket on $S$, equation (\ref{Poisson}), we
obtain
\begin{equation}
(\varphi_{-t}^{X_{h}})^{\ast}\{(\varphi_{t}^{X_{h}})f_{1},(\varphi_{t}^{X_{h}%
})f_{2}\}=\{f_{1},f_{2}\}\text{ for all }f_{1},f_{2},h\in C^{\infty}(S)
\label{invariance0}%
\end{equation}

Let $\mathcal{H}(S)$ denote the family of all Hamiltonian vector fields on $S
$.

\begin{description}
\item [Proposition 4.]The family $\mathcal{H}(S)$ is locally complete.
\end{description}%

\noindent
\textbf{Proof}. For $X_{f}\in\mathcal{H}(S)$, let $\varphi_{t}$ denote the
local one-parameter group of local diffeomorphisms of $S$ generated by $X$.
Suppose $U$ is the domain of $\varphi_{t}$ and $V$ is its range. In other
words $\varphi_{t}$ maps $U\,$ diffeomorphically onto $V$. For each $X_{h}%
\in\mathcal{H}(S)$, $\varphi_{t\ast}X_{h}$ is in $\mathrm{Der}(C^{\infty
}(V)).$ We have shown in the proof of Theorem 3 that $\varphi_{t\ast}X_{h}$ is
a vector field on $V$.

It follows from equation (\ref{push-forward}) that, for each $k\in C^{\infty
}(S)$,
\[
(\varphi_{t\ast}X_{h})\cdot(k\mid V)=\varphi_{-t}^{\ast}(X_{h}\cdot
(\varphi_{t}^{\ast}k))=\varphi_{-t}^{\ast}\{h,\varphi_{t}^{\ast}%
k\}=\{\varphi_{-t}^{\ast}h,k\mid V\},
\]
where we have taken into account equation (\ref{invariance0}). Hence,
$(\varphi_{t\ast}X_{h})$ is the inner derivation of $C^{\infty}(V)$
corresponding to the restriction of $\varphi_{-t}^{\ast}h=h\raisebox
{2pt}{$\scriptstyle\circ\, $}\varphi_{-t}$ to $V$.

For every $x\in V$, there exists an open neighbourhood $W$ of $x$ such that
$\bar{W}\subset V$. Let $\tilde{h}\in C^{\infty}(S)$ be such that $\tilde
{h}\mid W=\varphi_{-t}^{\ast}h\mid W.$ It follows that $(\varphi_{t\ast}%
X_{h})\mid W=X_{\tilde{h}}\mid W$. Hence, $\mathcal{H}(S)$ is complete.
\hfill$\Box$

Proposition 4 and Theorem 3 imply that orbits of the family $\mathcal{H}(S)$
of Hamiltonian vector fields on $S$ give rise to a singular foliation of $S$.
Moreover, local flows of Hamiltonian vector fields of $G$-invariant functions
on $P$ preserve local manifolds $P_{K}$, \cite{cushman-sniatycki}. Hence,
every Hamiltonian vector field $X_{h}$ on $S$, corresponding to $h\in
C^{\infty}(S)$, is strongly stratified with respect to the stratification of
$S$ by orbit type. Therefore, orbits of $\mathcal{H}(S)$ are contained in
strata of the stratification of $S$ by orbit type. Moreover, each orbit of
$\mathcal{H}(S)$ is a symplectic manifold, (\cite{liebermann-marle} p. 130).

\section{Subcartesian Poisson spaces}

In this section we generalize the notion of a Poisson manifold to a
subcartesian space.

Let $S$ be a subcartesian space. It will be called a Poisson space if its
differential structure $C^{\infty}(S)$ has a Poisson algebra structure and
inner derivations are vector fields on $S$. We denote by $\{f_{1},f_{2}\}$ the
Poisson bracket of $f_{1}$ and $f_{2}$ in $C^{\infty}(S),$ and assume that the
map
\[
C^{\infty}(S)\times C^{\infty}(S)\rightarrow C^{\infty}(S):(f_{1}%
,f_{2})\mapsto\{f_{1},f_{2}\}
\]
is bilinear, antisymmetric and satisfies both the Jacobi identity
\begin{equation}
\{f_{1},\{f_{2},f_{3}\}\}+\{f_{2},\{f_{3},f_{1}\}\}+\{f_{3},\{f_{1}%
,f_{2}\}\text{ for all }f_{1},f_{2},f_{3}\in C^{\infty}(S),
\label{Jacobi identity}%
\end{equation}
and the derivation condition
\begin{equation}
\{f_{1},f_{2}f_{3}\}=\{f_{1},f_{2}\}f_{3}+\{f_{1},f_{3}\}f_{2}\text{ for all
}f_{1},f_{2},f_{3}\in C^{\infty}(S). \label{derivation1}%
\end{equation}
Let
\begin{equation}
\mathcal{H}(S)=\left\{  X_{f}:C^{\infty}(S)\rightarrow C^{\infty
}(S):h\rightarrow X_{f}\cdot h=\{f,h\}\right\}  . \label{I}%
\end{equation}
It follows from equation (\ref{derivation1}) that each $X_{f}\in
\mathcal{H}(S)$ is a derivation of $C^{\infty}(S).$ Since derivations on a
subcartesian space needs not be vector fields, we make an additional
assumption that, for each $f\in C^{\infty}(S),$ the derivation $X_{f}%
\in\mathcal{H}(S)$ is a vector field on $S$. The vector field $X_{f}$ is
called the Hamiltonian vector field of $f$. The Jacobi identity implies that
$\mathcal{H}(S)$ is a Lie algebra subalgebra of $\mathrm{Der}C^{\infty}(S)$.
The Lie bracket on $\mathcal{H}(S)$ satisfies the identity
\[
\lbrack X_{f_{1}},X_{f_{2}}]=X_{\{f_{1},f_{2}\}}\text{ for all }f_{1},f_{2}\in
C^{\infty}(S)\text{. }%
\]
We shall refer to $\mathcal{H}(S)$ as the family of Hamiltonian vector
fields\textit{\ }on $S.$

\begin{description}
\item [Lemma 12.]For each open subset $U$ of a Poisson space $S$, the Poisson
bracket on $C^{\infty}(S)$ induces a Poisson bracket on $C^{\infty}(U).$
\end{description}

\noindent\textbf{Proof}. Let $h_{1},h_{2}\in C^{\infty}(U).$ Since $U$ is open
in $S$, for each point $x\in U$, there exists a neighbourhood $U_{x}$ of $x$
in $U$ and functions $f_{1,x},$ $f_{2,x}\in C^{\infty}(S)$ such that
$h_{1}\mid U_{x}=f_{1,x}\mid U_{x}$ and $h_{2}\mid U_{x}=f_{2,x}\mid U_{x}$.
We define
\begin{equation}
\{h_{1},h_{2}\}(x)=\{f_{1,x},f_{2,x}\}(x), \label{restricted}%
\end{equation}
where the right hand side is the value at $x\in U_{x}\subseteq S$ of the
Poisson bracket of functions in $C^{\infty}(S)$. We have to show that the
right hand side of equation (\ref{restricted}) is independent of the choice of
$U_{x}$ and $f_{1,x}$ and $f_{2,x}.$ Let $U_{x}^{\prime}$ be another open
neighbourhood of $x$ in $U$ and $f_{1,x}^{\prime},$ $f_{2,x}^{\prime}\in
C^{\infty}(S)$ such that $h_{1}\mid U_{x}^{\prime}=f_{1,x}^{\prime}\mid
U_{x}^{\prime}$ and $h_{2}\mid U_{x}^{\prime}=f_{2,x}^{\prime}\mid
U_{x}^{\prime}$. Then,
\begin{equation}
f_{1,x}\mid U_{x}\cap U_{x}^{\prime}=f_{1,x}^{\prime}\mid U_{x}\cap
U_{x}^{\prime}\text{ and }f_{2,x}\mid U_{x}\cap U_{x}^{\prime}=f_{2,x}%
^{\prime}\mid U_{x}\cap U_{x}^{\prime}. \label{comparison}%
\end{equation}
Hence, $k_{1,x}=f_{1,x}^{\prime}-f_{1,x}$ and $k_{2,x}=f_{2,x}^{\prime
}-f_{2,x}$ vanish on $U_{x}\cap U_{x}^{\prime}.$ Moreover,
\[
\{f_{1,x}^{\prime},f_{2,x}^{\prime}\}(x)=\{f_{1,x}+k_{1,x},f_{2,x}%
+k_{2,x}\}(x)=\{f_{1,x},f_{2,x}\}(x)
\]
because $\{f_{1,x},k_{2,x}\}$, $\{f_{2,x},k_{1,x}\}$ and $\{k_{1,x},k_{2,x}\}
$ vanish on $U_{x}\cap U_{x^{\prime}}.$ This proves that the Poisson bracket
on $C^{\infty}(U)$ is well defined by equation (\ref{restricted}). Moreover,
it is bilinear, antisymmetric and satisfies equations (\ref{Jacobi identity})
and (\ref{derivation1}) because the Poisson bracket on $C^{\infty}(S)$ has
these properties. \hfill$\Box$

\begin{description}
\item [Lemma 13.]The Poisson bracket $\{f_{1},f_{2}\}$ on $C^{\infty}(S)$ is
invariant under the local one-parameter groups of local diffeomorphisms of $S$
generated by Hamiltonian vector fields.
\end{description}

\noindent\textbf{Proof}. For $X_{f}\in\mathcal{H}(S)$, let $\varphi_{t}$
denote the local one-parameter group of local diffeomorphisms of $S$ generated
by $X$. Suppose $U$ is the domain of $\varphi_{t}$ and $V$ is its range. In
other words $\varphi_{t}$ maps $U\,$ diffeomorphically onto $V.$ We want to
show that
\begin{equation}
\varphi_{-t}^{\ast}\{\varphi_{t}^{\ast}f_{1},\varphi_{t}^{\ast}f_{2}%
\}=\{f_{1},f_{2}\}\text{ for all }f_{1},f_{2}\in C^{\infty}(S).
\label{invariance}%
\end{equation}

Differentiating the left-hand side with respect to $t$ we get
\begin{align*}
&  \frac{d}{dt}\varphi_{-t}^{\ast}\{\varphi_{t}^{\ast}f_{1},\varphi_{t}^{\ast
}f_{2}\}=-\varphi_{-t}^{\ast}\left(  X\cdot\{\varphi_{t}^{\ast}f_{1}%
,\varphi_{t}^{\ast}f_{2}\}\right)  +\\
&  +\varphi_{-t}^{\ast}\{X\cdot(\varphi_{t}^{\ast}f_{1}),\varphi_{t}^{\ast
}f_{2}\}+\varphi_{-t}^{\ast}\{\varphi_{t}^{\ast}f_{1},X\cdot(\varphi_{t}%
^{\ast}f_{2})\}\\
&  =\varphi_{-t}^{\ast}\left(  -\{\{f,\{\varphi_{t}^{\ast}f_{1},\varphi
_{t}^{\ast}f_{2}\}\}+\{\{f,\varphi_{t}^{\ast}f_{1}\},\varphi_{t}^{\ast}%
f_{2}\}+\{\varphi_{t}^{\ast}f_{1},\{f,\varphi_{t}^{\ast}f_{2}\}\}\right)  =0
\end{align*}
because of Jacobi identity (\ref{Jacobi identity}). Since $\varphi_{0}^{\ast
}f_{i}=f_{i}$, it follows that
\[
\varphi_{-t}^{\ast}\{\varphi_{t}^{\ast}f_{1},\varphi_{t}^{\ast}f_{2}%
\}=\varphi_{-0}^{\ast}\{\varphi_{0}^{\ast}f_{1},\varphi_{0}^{\ast}%
f_{2}\}=\{f_{1},f_{2}\},
\]
which completes the proof. \hfill$\Box$

\begin{description}
\item [Lemma 14.]Let $S$ be a subcartesian Poisson space. The family
$\mathcal{H}(S)$ of Hamiltonian vector fields is locally\ complete.
\end{description}

\noindent Proof of this lemma is identical to the proof of Proposition 4.

\begin{description}
\item [Theorem 9.]Let $S$ be a subcartesian Poisson space. Orbits of the
family $\mathcal{X}(S)$ of all vector fields on $S$ are Poisson manifolds.
Orbits of the family $\mathcal{H}(S)$ of Hamiltonian vector fields are
symplectic manifolds. For each Poisson leaf $M$, orbits of $\mathcal{H}(S)$
contained in $M$ give rise to a singular foliation of $M$ by symplectic leaves.
\end{description}

\noindent\textbf{Proof}. Let $M$ be an orbit of $\mathcal{X}(S)$. By Theorem
4, it is a smooth manifold. Let $C^{\infty}(M)$ be the space of smooth
functions on $M$ defined in terms of the manifold structure of $M$. Let
$\mathcal{T}$ denote the manifold topology of $M$ described in section 6. A
function $h:M\rightarrow\mathbb{R}$ is in $C^{\infty}(M)$ if and only if, for
each $x\in M$, there exists $V\in\mathcal{T}$ such that $x\in V$ and there
exists a function $f_{x}\in C^{\infty}(S)$ such that $h\mid V=f_{x}\mid V.$

For $h_{1},h_{2}\in C^{\infty}(M)$ and $x\in M$, let $V\in\mathcal{T},$ and
$f_{1,x},$ $f_{2,x}\in C^{\infty}(S)$ be such that $h_{i}\mid V=f_{i,x}\mid V$
for $i=1,2.$ We define $\{h_{1},h_{2}\}_{M}$ by the requirement that
\begin{equation}
\{h_{1},h_{2}\}_{M}\mid V=\{f_{1,x},f_{2,x}\}\mid V. \label{bracket}%
\end{equation}
We have to show that the right-hand side of equation (\ref{bracket}) is
independent of the choice of $f_{1,x},$ and $f_{2,x}$ in $C^{\infty}(S)$.
Suppose that $f_{1,x}^{\prime},$ $f_{2,x}^{\prime}\in C^{\infty}(S)$ satisfy
the condition $h_{i}\mid V=f_{i,x}^{\prime}\mid V$ for $i=1,2.$ Let
$k_{i,x}=(f_{i,x}^{\prime}-f_{i,x}).$ Then, $k_{i,x}\mid V=0,$ and
\begin{align*}
\{f_{1,x}^{\prime},f_{2,x}^{\prime}\}  &  \mid V=\{f_{1,x}+k_{1,x}%
,f_{2,x}+k_{2,x}\}\mid V\\
&  =\left(  \{f_{1,x},f_{2,x}\}+\{f_{1,x},k_{2,x}\}+\{k_{1,x},f_{2,x}%
\}+\{k_{1,x},k_{2,x}\}\right)  \mid V\\
&  =(\{f_{1,x},f_{2,x}\}+(X_{f_{1,x}}\cdot k_{2,x})-(X_{f_{2,x}}\cdot
k_{1,x})+(X_{k_{1,x}}\cdot k_{2,x}))\mid V.
\end{align*}
Since $M$ is an orbit of $\mathcal{X}(S),$ $V$ is an open subset of $M$, and
$k_{1,x}\mid V=0$ and $k_{2,x}\mid V=0,$ it follows that $(X\cdot k_{1,x})\mid
V=(X\cdot k_{2,x})\mid V=0$ for all $X\in\mathcal{X}(S)$. Hence,
$\{f_{1,x}^{\prime},f_{2,x}^{\prime}\}\mid V=\{f_{1,x},f_{2,x}\}\mid V=0$ and
$\{h_{1},h_{2}\}_{M}$ is well defined.

It follows from equation (\ref{bracket}) that $\{h_{1},h_{2}\}_{M}\in
C^{\infty}(M).$ The Poisson bracket properties of $\{.,.\}_{M}$ follow from
the corresponding properties of the Poisson bracket on $C^{\infty}(S).$

By Lemma 14, the family $\mathcal{H}(S)$ of Hamiltonian vector fields on $S$
is complete. Hence, its orbits give rise to a singular foliation of $S$.
Theorem 4 implies that each orbit of $\mathcal{H}(S)$ is contained in an orbit
of $\mathcal{X}(S)$.

We have shown that each orbit $M$ of $\mathcal{X}(S)$ is a Poisson manifold.
Orbits of $\mathcal{H}(S)$ contained in $M$ coincide with orbits of the family
of Hamiltonian vector fields on $M$, which gives rise to a foliation of $M$ by
symplectic leaves of $M,$ (\cite{liebermann-marle} p. 130). \hfill$\Box$

\bigskip

Let $S$ be a Poisson space, and $G$ be a connected Lie group with Lie algebra
$\frak{g}$. We denote by
\[
\Phi:G\times S\rightarrow S:(g,x)\mapsto\Phi(g,x)\equiv\Phi_{g}(x)=gx
\]
an action of $G$ on $S$. We assume that this action is smooth, which implies
that, for each $g\in G$, the map $\Phi_{g}:S\rightarrow S$ is a
diffeomorphism. Moreover, we assume that, for every $g\in G,$
\[
\Phi_{g}^{\ast}:C^{\infty}(S)\rightarrow C^{\infty}(S):f\mapsto\Phi_{g}^{\ast
}(f)=f\raisebox{2pt}{$\scriptstyle\circ\, $}\Phi_{g}%
\]
is an automorphism of the Poisson algebra structure of $C^{\infty}(S)$. In
other words,
\[
\Phi_{g}^{\ast}\{f_{1},f_{2}\}=\{\Phi_{g}^{\ast}f_{1},\Phi_{g}^{\ast}%
f_{2}\}\text{ for all }g\in G\text{ and }f_{1},f_{2}\in C^{\infty}(S).
\]
Finally, we assume that the action $\Phi$ is proper.

For each $\xi\in\frak{g}$, we denote by $X^{\xi}$ the vector field on $S$
generating the action on $S$ of the one-parameter subgroup $\exp t\xi$ of $G
$. Clearly, $X^{\xi}\in X(S)$ for all $\xi\in\frak{g}$. Since $G$ is
connected, its action on $S$ is generated by the action of all one-parameter
subgroups. Hence, each Poisson manifold of $S$ is invariant under the action
of $G$ on $S.$ We denote by $\Phi^{M}:G\times M\rightarrow M$ the induced
action of $G$ on a Poisson manifold $M$. The assumptions on the action of $G$
on $S$ imply that the action $\Phi^{M}$ is smooth and proper. Moreover, it
preserves the Poisson algebra structure of $C^{\infty}(M).$

Let $R_{M}=M/G$ denote the space of $G$-orbits in $M$ and $\rho_{M}%
:M\rightarrow R_{M}$ the orbit map. Since $M$ is a manifold, and the action of
$G$ on $M$ is proper, it follows that $R_{M}$ is a stratified space which can
be covered by open sets, each of which is diffeomorphic to an open subset of a
semi-algebraic set, (\cite{cushman-sniatycki}, p. 727). Hence, $R_{M}$ is a
subcartesian space.

The differential structure $C^{\infty}(R_{M})$ is isomorphic to the ring
$C^{\infty}(M)^{G}$ of $G$-invariant smooth functions on $M$. Since the action
of $G$ on $M$ preserves the Poisson bracket $\{.,.\}_{M}$ on $C^{\infty}(M)$,
it follows that $C^{\infty}(M)^{G}$ is a Poisson subalgebra of $C^{\infty}%
(M)$. Hence, $C^{\infty}(R_{M})$ inherits a Poisson bracket $\{.,.\}_{R_{M}}$.

For $h\in C^{\infty}(R_{M})$, let $f=\rho_{M}^{\ast}h\in C^{\infty}(M)^{G}.$
Let $X_{h}$ be the derivation of $C^{\infty}(R_{M})$ given by $X_{h}\cdot
h^{\prime}=\{h,h^{\prime}\}_{R_{M}}$ for all $h^{\prime}\in C^{\infty}(R_{M}%
)$. Similarly, $X_{f}$ is the derivation of $C^{\infty}(M)$ given by
$X_{f}\cdot f^{\prime}=\{f,f^{\prime}\}_{M}$ for all $f^{\prime}\in C^{\infty
}(M)$. The vector field $X_{f}$ generates a one-parameter group $\varphi_{t}$
of local diffeomorphisms of $M$ which commutes with the action of $G$ on $M$.
Hence, $\varphi_{t}$ induces a local one-parameter group of local
diffeomorphisms $\psi_{t}$ of the orbit space $R_{M}=M/G$ such that $\psi
_{t}\raisebox{2pt}{$\scriptstyle\circ\, $}\rho_{M}=\rho_{M}\raisebox
{2pt}{$\scriptstyle\circ\, $}\varphi_{t}$. The local group $\psi_{t}$ is
generated by $X_{h}$. Hence, $X_{h}$ is a vector field on $R_{M}$.

It follows from the above discussion that the orbit space $R_{M}$ is a
subcartesian Poisson space. Hence, we can apply the results of Theorem 9.

\begin{description}
\item [Proposition 5.]Let $M$ be a Poisson manifold, and $R_{M}=M/G$ be the
orbit space of a properly acting Lie symmetry group $G$ of the Poisson
structure on $M.$ Then $R_{M}$ is a subcartesian Poisson space. It is a
stratified space. Strata of $R_{M}$ are orbits of the family $\mathcal{X}%
(R_{M})$ of all vector fields on $R_{M}$. Each stratum is a Poisson manifold.
The singular foliation of $R_{M}$ by orbits of the Lie algebra $\mathcal{H}%
(R_{M})$ of Hamiltonian vector fields of $C^{\infty}(R_{M})$ gives rise to a
refinement of the stratification of $R_{M}$ by symplectic manifolds.
\end{description}

\noindent\textbf{Proof}. It follows from Theorem 9 that orbits of the family
$\mathcal{X}(R_{M})$ of all vector fields on $R_{M}$ are Poisson manifolds.
Stratification structure of $R_{M}$ and its smootly local triviality are
consequences of the properness of the action of $G$ on $M,$ \cite{D-K},
\cite{bierstone1}, \cite{bierstone2}. It follows from Theorem 7 and Theorem 9
that strata of $R_{M}$ are Poisson manifolds. Moreover, Theorem 9 implies that
orbits of the family $\mathcal{H}(R_{M})$ of Hamiltonian vector fields are
symplectic manifolds. The restriction of the singular foliation of $R_{M}$ by
symplectic manifolds to each stratum of $R_{M}$\ gives rise to a singular
foliation of the stratum by symplectic manifolds. \hfill$\Box$

\bigskip

Let $R=S/G$ be the space of $G$-orbits in $S$ and $\rho:S\rightarrow R$ the
orbit map. It is a differential space with differential structure $C^{\infty
}(R)$ isomorphic to the ring $C^{\infty}(S)^{G}$ of $G$-invariant smooth
functions on $S$. The Poisson algebra structure of $C^{\infty}(S)$ induces a
Poisson structure on $C^{\infty}(R).$ It follows from Corollary 3 and the
discussion preceding it, that $R$ is singularly foliated by Poisson manifolds,
and each Poisson leaf is singularly foliated by symplectic leaves. We do not
know if $R$ is a subcartesian space. Hence, we cannot assert that Poisson
leaves of $R$ are orbits of the family $\mathcal{X}(R)$ of all vector fields
on $R,$ or that symplectic leaves of $R$ are orbits of the family
$\mathcal{H}(R)$ of Hamiltonian vector fields $.$

\section{Almost complex structures}

In this section, we discuss almost complex structures defined on complete
families of vector fields on subcartesian spaces. We assume here that the
subcartesian space under consideration is paracompact. By a theorem of
Marshall, this assumption ensures the existence of partitions of unity,
\cite{marshall}.

Let $\mathcal{F}=\{X^{\alpha}\}_{\alpha\in A}$ be a complete family of vector
fields on a paracompact subcartesian space $S$. We denote by $\mathrm{Der}%
_{\mathcal{F}}(C^{\infty}(S))$ the submodule of derivations of $C^{\infty}(S)$
consisting of locally finite sums $\Sigma_{\alpha}f_{\alpha}X^{\alpha}$, where
$f_{\alpha}\in C^{\infty}(S)$, $X^{\alpha}\in\mathcal{F}$ and, for every $x\in
S$, there is an open neighbourhood $U$ of $x$ in $S$ such that $f_{\alpha
}X^{\alpha}\mid U=0$ for almost all $\alpha.$ Abusing somewhat the common
terminology, we shall refer to $\mathrm{Der}_{\mathcal{F}}(C^{\infty}(S))$ as
the module of derivations generated by $\mathcal{F}$.

\begin{description}
\item [Proposition 6.]For every complete family $\mathcal{F}$ of vector fields
on a paracompact subcartesian space $S$, the module $\mathrm{Der}%
_{\mathcal{F}}(C^{\infty}(S))$\ of derivations generated by $\mathcal{F}$ is a
Lie subalgebra of the Lie algebra of all derivations of $C^{\infty}(S).$
\end{description}%

\noindent
\textbf{Proof}. Recall that the completeness of $\mathcal{F}=\{X^{\alpha
}\}_{\alpha\in A}$ implies that, for every $\alpha,$ $\beta$, $t,$ and $x$ for
which $\varphi_{t\ast}^{\alpha}X^{\beta}(x)$ is defined, there exists an open
neighbourhood $U$ of $x$ and $\gamma\in A$ such that $\varphi_{t\ast}^{\alpha
}X^{\beta}\mid U=X^{\gamma}\mid U.$ Hence, there exists $Z_{U}^{\alpha\beta
}\in\mathrm{Der}_{\mathcal{F}}(C^{\infty}(S))$ such that $[X^{\alpha}%
,X^{\beta}]\mid U=Z_{U}^{\alpha\beta}\mid U$. In this way we get an open cover
$\mathcal{U=}\{U\}$ of $S$. By shrinking open sets $U$, if necessary, we may
assume that the covering $\mathcal{U=}\{U\}$ is locally finite. Since $S$ is
paracompact, there exists a partition of unity $\{f_{U}\}_{U\in\mathcal{U}}$
subordinate to this covering, \cite{marshall}. Hence, $[X^{\alpha},X^{\beta
}]\mid U=\Sigma_{U}f_{U}(Z_{U}^{\alpha\beta}\mid U)$, where the sum on the
right-hand side is locally finite. This implies that $[X^{\alpha},X^{\beta
}]\in\mathrm{Der}_{\mathcal{F}}(C^{\infty}(S))$.

If $X=\Sigma_{\alpha}f_{\alpha}X^{\alpha}$ and $Y=\Sigma_{\beta}h_{\beta
}X^{\beta}$ are in $\mathrm{Der}_{\mathcal{F}}(C^{\infty}(S)),$ then
\begin{align*}
\lbrack X,Y]  &  =\left[
{\textstyle\sum}
{}_{\alpha}f_{\alpha}X^{\alpha},%
{\textstyle\sum}
{}_{\beta}h_{\beta}X^{\beta}\right] \\
&  =%
{\textstyle\sum}
{}_{\alpha\beta}\left(  f_{\alpha}h_{\beta}[X^{\alpha},X^{\beta}]+f_{\alpha
}(X^{\alpha}\cdot h_{\beta}X^{\beta}-h_{\beta}(X^{\beta}\cdot f_{\alpha
})X^{\alpha}\right)  .
\end{align*}
Since the sum is locally finite, it implies that $[X,Y]\in\mathrm{Der}%
_{\mathcal{F}}(C^{\infty}(S)).$ \hfill$\Box$

An\textit{\ }almost complex structure on a complete family $\mathcal{F}$ of
vector fields on $S$ is a $C^{\infty}(S)$ module automorphism $J:\mathrm{Der}%
_{\mathcal{F}}(C^{\infty}(S))\rightarrow\mathrm{Der}_{\mathcal{F}}(C^{\infty
}(S))$ such that $J^{2}=-1.$ Since $J:\mathrm{Der}_{\mathcal{F}}(C^{\infty
}(S))\rightarrow\mathrm{Der}_{\mathcal{F}}(C^{\infty}(S))$ is a $C^{\infty
}(S)$ module automorphism, it implies that, for each orbit $M$ of
$\mathcal{F}$, it gives rise to a linear map $J_{M}:TM\rightarrow TM$.
Moreover, $J^{2}=-1$ implies that $J_{M}^{2}=-1.$ Hence, an almost complex
structure on a complete family $\mathcal{F}$ of vector fields on $S$ induces
an almost complex structure on each orbit of $\mathcal{F}$.

Since $\mathrm{Der}_{\mathcal{F}}(C^{\infty}(S))$ is a Lie algebra, we may
consider the torsion $N$ of $J$ defined as follows. For $X,Y\in\mathrm{Der}%
_{\mathcal{F}}(C^{\infty}(S))$, let
\begin{equation}
N(X,Y)=2\left\{  [JX,JY]-J[JX,Y]-J[X,JY]-[X,Y]\right\}  . \label{N}%
\end{equation}

\begin{description}
\item [Lemma 15.]The torsion of $J$ is a skew symmetric bilinear mapping
$N:\mathrm{Der}_{\mathcal{F}}(C^{\infty}(S))\otimes\mathrm{Der}_{\mathcal{F}%
}(C^{\infty}(S))\rightarrow\mathrm{Der}_{\mathcal{F}}(C^{\infty}(S))$ such
that $N(fX,hY)=fhN(X,Y)$ for all $X,Y$ in $\mathrm{Der}_{\mathcal{F}%
}(C^{\infty}(S))$ and $f,h\in C^{\infty}(S).$
\end{description}%

\noindent
\textbf{Proof}. Skew symmetry and bilinearity of $N$ are self-evident. For
every $X,Y\in\mathrm{Der}_{\mathcal{F}}(C^{\infty}(S))$ and $f,h\in C^{\infty
}(S).$
\begin{align*}
\lbrack JfX,JhY]  &  =fh[JX,JY]+f((JX)\cdot h)JY-h((JY)\cdot f)JX\\
J[fX,JhY]  &  =fhJ[X,JY]-f(X\cdot h)Y-h((JY)\cdot f)JX.\\
J[JfX,hY]  &  =fhJ[JX,Y]+f((JX)\cdot h)JY+h(Y\cdot f)X\\
\lbrack fX,hY]  &  =fh[X,Y]+f(X\cdot h)Y-h(Y\cdot f)X.
\end{align*}
Hence, $N(fX,hY)=fhN(X,Y),$ which completes the proof. \hfill$\Box$

\bigskip

The almost complex structure $J$ has eigenvalues $\pm i$ because $J^{2}=-1.$
Eigenspaces of $J$ are contained in the complexification $\mathrm{Der}%
_{\mathcal{F}}(C^{\infty}(S))\otimes\mathbb{C}$ of $\mathrm{Der}_{\mathcal{F}%
}(C^{\infty}(S)).$ For every $X\in\mathrm{Der}_{\mathcal{F}}(C^{\infty}(S))$,
\begin{align*}
J(X-iJX)  &  =JX-iJ^{2}X=JX+iX=i(X-iJX),\\
J(X+iJX)  &  =JX+iJ^{2}X=JX-iX=-i(X+iJX).
\end{align*}
Hence,
\[
\mathrm{Der}_{\mathcal{F}}(C^{\infty}(S))^{\pm}=\left\{  X-(\pm i)JX\mid
X\in\mathrm{Der}_{\mathcal{F}}(C^{\infty}(S))\right\}
\]
are eigenspaces of $J$ corresponding to eigenvalues $\pm i$, respectively.

\begin{description}
\item [Lemma 16.]Eigenspaces $\mathrm{Der}_{\mathcal{F}}(C^{\infty}(S))^{\pm}$
of $J$ are closed under the Lie bracket if and only if the torsion $N$ of $J$ vanishes.
\end{description}%

\noindent
\textbf{Proof}. For every $X$ and $Y$ in $\mathrm{Der}_{\mathcal{F}}%
(C^{\infty}(S)),$ we have
\begin{align*}
&  [X-(\pm i)JX,Y-(\pm i)JY]\\
&  =[X,Y]-(\pm i)[JX,Y]-(\pm i)[X,JY]+(\pm i)^{2}[JX,JY]\\
&  =[X,Y]-(\pm i)[JX,Y]-(\pm i)[X,JY]-[JX,JY]\\
&  =-N(X,Y)-J[JX,Y]-J[X,JY]-(\pm i)[JX,Y]-(\pm i)[X,JY]\\
&  =-N(X,Y)-(\pm i)\left\{  \left(  [JX,Y]+[X,JY]\right)  -(\pm iJ)\left(
[JX,Y]+[X,JY]\right)  \right\}  .
\end{align*}
Hence, $[X-(\pm i)JX,Y-(\pm i)JY]\in\mathrm{Der}_{\mathcal{F}}(C^{\infty
}(S))^{\pm}$ if and only if $N(X,Y)=0.$ \hfill$\Box$

\begin{description}
\item [Theorem 10.]Let $J$ be an almost complex structure on a complete family
$\mathcal{F}$ of vector fields on a paracompact subcartesian space $S$. Every
orbit $M$ of $\mathcal{F}$ admits a complex analytic structure such that
$\mathrm{Der}_{\mathcal{F}}(C^{\infty}(S))^{+}\mid M$ spans the distribution
of holomorphic directions on $TM\otimes\mathbb{C}$ if and only if the torsion
$N$ of $J$ vanishes$.$
\end{description}%

\noindent
\textbf{Proof}. For each orbit $M$ of $\mathcal{F}$, the restriction of $N$ to
$M$ is the torsion tensor $N_{M}$ of the almost complex structure $J_{M}$ on
$M$. Suppose that $M$ admits a complex analytic structure such that
$\mathrm{Der}_{\mathcal{F}}(C^{\infty}(S))^{+}\mid M$ spans the distribution
of holomorphic directions on $TM\otimes\mathbb{C}$. For a complex manifold
$M$, the distribution of holomorphic directions on $TM\otimes\mathbb{C}$ is
involutive. Hence, $\mathrm{Der}_{\mathcal{F}}(C^{\infty}(S))^{+}\mid M$ is
closed under the Lie bracket of vector fields on $TM\otimes\mathbb{C}$. This
implies that $N_{M}=0.$ Therefore, if every orbit $M$ of $\mathcal{F}$ admits
a complex analytic structure such that $\mathrm{Der}_{\mathcal{F}}(C^{\infty
}(S))^{+}\mid M$ spans the distribution of holomorphic directions on
$TM\otimes\mathbb{C}$, then the torsion $N$ vanishes.

Suppose now that $N=0$. Then $N_{M}=0$ for every orbit $M$ of $\mathcal{F}$.
It follows that the almost complex structure $J_{M}$ on $M$ is integrable. By
the Newlander-Nirenberg theorem, there exists a complex analytic structure on
$M$ such that $\mathrm{Der}_{\mathcal{F}}(C^{\infty}(S))^{+}\mid M$ spans the
distribution of holomorphic directions on $TM\otimes\mathbb{C}$,
\cite{newlander-nirenberg}. Naively speaking, one could say that this complex
structure is obtained by patching local diffeomorphisms from $\mathbb{C}^{n}$
to $M$ in a manner similar to that used in the proof of Theorem 3.\hfill$\Box$

Assuming that the torsion tensor $N$ of $J$ vanishes, we are going to
cha\-racterize smooth functions on $S$ which are holomorphic on each orbit
of\thinspace$\mathcal{F}$. Let $C^{\infty}(S)^{\mathbb{C}}=C^{\infty
}(S)\otimes\mathbb{C}$ be the complexification of $C^{\infty}(S).$ Each
function in $C^{\infty}(S)^{\mathbb{C}}$ is of the form $f+ih$, where $f,h\in
C^{\infty}(S)$. Such a function is holomorphic on\thinspace each orbit
of\thinspace$\mathcal{F}$ if it is annihilated by derivations in
$\mathrm{Der}_{\mathcal{F}}(C^{\infty}(S))^{-}.$ In other words, a function
$f+ih$ is holomorphic on orbits of $\mathcal{F}$ if it satisfies the
differential equation $(X+iJX)(f+ih)=0$ for each $X\in\mathcal{F}.$ Separating
real and imaginary parts we get a singularly foliated version of
Cauchy-Riemann equations
\[
X\cdot f-(JX)\cdot h=0\text{ and }(JX)\cdot f+X\cdot h=0\text{ for all }%
X\in\mathcal{F}.
\]
It should be noted that these equations may have very few solutions which are
globally defined. This is why one usually employs sheaves in the study of
holomorphic functions, \cite{spallek1}.

Combining the results of the last two sections, we can describe subcartesian
Poisson-K\"{a}hler spaces.\bigskip

\noindent\textbf{Example 3}. Let $S$ be a subcartesian Poisson space. It
follows from the definition of Hamiltonian vector fields on $S$ that the
Poisson bracket $\{f,h\}$ on $C^{\infty}(S)$ satisfies the relations
\[
X_{f}\cdot h=\{f,h\}=-\{h,f\}=-X_{h}\cdot f.
\]
Hence, there exists a skew symmetric form $\Omega$ on $H(S)$ with values in
$C^{\infty}(S)$ such that
\[
\Omega(X_{f},X_{h})=\{f,h\}\text{ }\forall\text{ }f,h\in C^{\infty}(S)\text{.}%
\]
By Theorem 9, every orbit $M$ of $\mathcal{H}(S)$ is a symplectic manifold.
The symplectic form $\Omega_{M}$ of $M$ is given by the restriction of
$\Omega$ to $M.$ Let $J:\mathrm{Der}_{\mathcal{H}(S)}C^{\infty}(S)\rightarrow
\mathrm{Der}_{\mathcal{H}(S)}C^{\infty}(S)$ be an almost complex structure on
$\mathcal{H}(X)$ such that,
\[
\Omega(JX,JY)=\Omega(Y,X)\text{ }\forall\text{ }X,Y\in\mathrm{Der}%
_{\mathcal{H}(S)}C^{\infty}(S),
\]
for all $X\in\mathrm{Der}_{\mathcal{H}(S)}C^{\infty}(S).$ Define
\[
g(X,Y)=\Omega(JX,Y)
\]
for all $X,Y\in\mathrm{Der}_{\mathcal{H}(S)}C^{\infty}(S)$. It is symmetric
because
\begin{align*}
g(Y,X)  &  =\Phi(JY,X)=\Phi(J^{2}Y,JX)=-\Phi(Y,JX)=\Phi(JX,Y)\\
&  =g(X,Y).
\end{align*}
If $g(X,Y)$ is positive definite, that is $g(X,X)(x)=0$ only if $X(x)=0$, then
its restriction to every orbit $M$ of $H(S)$ defines a Riemannian metric
$g_{M}$ on $M$. For every $X,Y\in H(S)$ and $x\in M$, $g_{M}(X(x),Y(x))=\Omega
_{M}(JX(x),Y(x))$ and the form $\Omega_{M}$ is closed. Hence, $M$ is an almost
K\"{a}hler manifold. If \ the torsion $N$ of $J$ vanishes, then every orbit
$M$ of $\mathcal{H}(S)$ is a K\"{a}hler manifold. This example is a
generalization to subcartesian spaces of K\"{a}hler stratified spaces studied
by Huebschmann, \cite{huebschmann}.\bigskip

\noindent{\Large Acknowledgments.\ }

\medskip

I would like to thank Larry Bates, Cathy Beveridge, Richard Cushman, Hanafi
Farahat, Johannes Huebschmann and James Montaldi for stimulating discussions,
and Mark Roberts for an invitation to the Warwick Symposium on Geometric
Mechanics and Symmetry at the Mathematics Research Centre, University of
Warwick, where I began working on this paper. I would also like to thank the
Referee whose comments and corrections have lead to an improvement of the paper.

\end{document}